\documentclass[12pt]{article}
\usepackage{amsmath}
\usepackage{amsfonts}
\usepackage{mathrsfs}
\usepackage{amssymb}
\usepackage[usenames]{color}
\usepackage[numbers,sort&compress]{natbib}
\usepackage{subfigure}
\usepackage[all]{xy}

\def\cl{\centerline}
\def\ni{\noindent}
\def\vs{\vspace*}

\def\a{\alpha}
\def\b{\beta}
\def\D{\Delta}
\def\L{\Lambda}

\def\CVir{{\frak {Vir}}}
\def\CB{{\frak {B}}(p)}

\def\BB{{\cal B}(q)}
\def\sp{{\rm{span}}}
\def\deg{{\rm deg}}

\def\NS{{\frak {NS}}}
\def\SB{{\frak {S}}(p)}

\def\K2{K_2}
\def\KB{{\frak {K}}(p)}
\def\KBB{{\frak {K}}}

\def\Z{\mathbb{Z}}

\def\C{\mathbb{C}}

\def\QED{\hfill$\Box$}

\numberwithin{equation}{section}
\newtheorem{theo}{Theorem}[section]
\newtheorem{defi}[theo]{Definition}
\newtheorem{coro}[theo]{Corollary}
\newtheorem{lemm}[theo]{Lemma}
\newtheorem{rema}[theo]{Remark}
\newtheorem{prop}[theo]{Proposition}

\newtheorem{exam}{Example}

\begin{document}
\begin{center}
{\bf\large Classification of finite irreducible conformal modules over $N=2$ Lie conformal superalgebras of Block type}
\footnote{
$^{*}$Corresponding author: chgxia@cumt.edu.cn (C.~Xia).}
\end{center}
\vspace{7pt}

\cl{Chunguang Xia$^{\,*}$}
\vspace{7pt}

\cl{\small School of Mathematics, China University of Mining and Technology, Xuzhou 221116, China}
\vspace{7pt}

\cl{\small Email: chgxia@cumt.edu.cn}
\vspace{7pt}

\footnotesize
\noindent{{\bf Abstract.}
We introduce the $N=2$ Lie conformal superalgebras $\KB$ of Block type, and classify their finite irreducible conformal modules
for any nonzero parameter $p$. 
In particular, we show that such a conformal module admits a nontrivial
extension of a finite conformal module $M$ over $\K2$ if $p=-1$ and $M$ has rank $(2+2)$, where $\K2$ is an $N=2$ conformal subalgebra of $\KB$.
As a byproduct, we obtain the classification of finite irreducible conformal modules over a series of finite Lie conformal superalgebras ${\frak k}(n)$ for $n\ge1$.
Composition factors of all the involved reducible conformal modules are also determined.
\vs{5pt}

\ni{\bf Key words:} finite conformal module; $N=2$ conformal superalgebra; composition factor

\ni{\it Mathematics Subject Classification (2010):} 17B10; 17B65; 17B68; 17B69.}

\small
\section{\large{Introduction}}
The present paper is the third in our series of papers on representation theory
of infinite Lie conformal superalgebras with Cartan type conformal subalgebras, the first two of which are \cite{SXY} and \cite{X}.

Lie conformal superalgebras  \cite{K1} encode the singular part of the operator
product expansion of chiral fields in conformal field theory.
During the last two decades, many advances have been made in the theory of finite Lie conformal superalgebras \cite{FK,FKR,CK,CL,BKLR,BKL4,BKL3,MZ,Z1,Z2}.
A complete classification of finite simple Lie conformal superalgebras can be found in \cite{FK},
which includes Cartan type and current type Lie conformal superalgebras.
The most physically important Cartan type Lie conformal superalgebras include the Virasoro conformal algebra $K_0$, the Neveu-Schwarz
conformal algebra $K_1$, and the $N=2$ conformal superalgebra $K_2$.
For the classification of finite irreducible conformal modules (FICMs)
over finite simple Lie conformal superalgebras, see \cite{CK,CL,BKLR,BKL4,BKL3, MZ}.

The theory of infinite Lie conformal superalgebras, however, is only in the early stages of development.
Let $p$ be a nonzero complex number. In \cite{SXY}, we introduced the infinite Lie conformal algebra $\CB$, which
contains a Virasoro conformal subalgebra $K_0$ (denoted $\CVir$ in \cite{SXY}) and has close relation with Lie algebras of Block type \cite{SXX1,SXX2,XZ}.
In \cite{X}, we constructed the Lie conformal superalgebra $\SB$ as the super analogue of $\CB$, which
contains a Neveu-Schwarz conformal subalgebra $K_1$ (denoted $\NS$ in \cite{X}). It is well-known that $K_0\subset K_1\subset K_2$.
In the present paper, we shall construct the Lie conformal superalgebra $\KB$ such
that the following embedding diagram is commutative:
$$
\xymatrix{
  K_0 \ar[d]_{} \ar[r]^{} & K_1 \ar[d]_{} \ar[r]^{} & K_2 \ar[d]^{} \\
  \CB \ar[r]^{} & \SB \ar[r]^{} & \KB }
$$
Naturally, we refer to $\KB$'s as {\it $N=2$ Lie conformal superalgebras of Block type}.

As one can see later, $\KB=\KB_{\bar{0}}\oplus\KB_{\bar{1}}$, where the even part $\KB_{\bar{0}}$ has $\C[\partial]$-basis
$\{L_i,\,J_i\,|\,i\in\Z_+\}$ and the odd part $\KB_{\bar{1}}$ has $\C[\partial]$-basis $\{G^{\pm}_i\,|\,i\in\Z_+\}$,
satisfying $\lambda$-brackets
\begin{eqnarray}
\label{brackets-LL}[L_i\, {}_\lambda \, L_j] &\!\!\!=\!\!\!& ((i+p)\partial+(i+j+2p)\lambda) L_{i+j}, \\[-2pt]
\label{brackets-LJ}[L_i\, {}_\lambda \, J_j] &\!\!\!=\!\!\!& ((i+p)\partial+(i+j+p)\lambda) J_{i+j}, \\[-2pt]
\label{brackets-LG}[L_i\, {}_\lambda \, G^{\pm}_j] &\!\!\!=\!\!\!& ((i+p)\partial+(i+j+\frac{3}{2}p)\lambda) G^{\pm}_{i+j}, \\[-2pt]
\label{brackets-JG}[J_i\, {}_\lambda \, G^{\pm}_j] &\!\!\!=\!\!\!& \pm G^{\pm}_{i+j}, \\[-2pt]
\label{brackets-GG}[G^{+}_i\, {}_\lambda \, G^{-}_j] &\!\!\!=\!\!\!& ((2i+p)\partial+2(i+j+p)\lambda) J_{i+j}+2 L_{i+j}.
\end{eqnarray}
Other $\lambda$-brackets are given by the skew-symmetry or vanish.
It is worth to highlight some interesting features of $\KB$.
Firstly, each $\KB$ contains an $N=2$ conformal subalgebra.
Setting $L=\frac{1}{p}L_0$, $J=J_0$, $G^{\pm}=\frac{1}{\sqrt{p}}G^{\pm}_0$,
one can check that the subalgebra
\begin{equation}\label{equ-K2}
K_2=\C[\partial]L\oplus\C[\partial]J\oplus\C[\partial]G^{+}\oplus\C[\partial]G^{-}
\end{equation}
of $\KB$ is exactly the $N=2$ conformal superalgebra \cite{FK} (see also Subsection~2.2).
Secondly, the subalgebra of $\KB$ with $\C[\partial]$-basis $\{L_i,\,G_i=\frac{1}{\sqrt{2}}(G^{+}_i+G^{-}_i)\,|\,i\in\Z_+\}$
is isomorphic to $\SB$ \cite{X}.
These two features suggest that the above embedding diagram is commutative.
Thirdly, there are embedding relations among $\KB$'s. For any integer $n\ge 1$, $\KB$ can be embedded into $\KBB(np)$
via $L_i\mapsto \frac{1}{n}L'_{ni}$, $J_i\mapsto J'_{ni}$, $G^{\pm}_i\mapsto \frac{1}{\sqrt{n}}{G^{\pm}_{ni}}'$.
Finally, $\KBB(-n)$ contains a series of finite Lie conformal superalgebras as quotient algebras (cf.~\eqref{equ-quotiont-special})
\begin{equation}\label{equ-sn}
{\frak k}(n)=\KBB(-n)/\KBB(-n)_{\langle n+1 \rangle}.
\end{equation}
The special cases ${\frak k}(1)$ and ${\frak k}(2)$ will be referred to as {\it $N=2$ Heisenberg conformal superalgebra}
and  {\it $N=2$ Schr$\ddot{\rm o}$dinger conformal superalgebra}, respectively. See Subsection~2.3 for more details.

Our main goal in this paper is to classify FICMs over $\KB$.
A complete classification of FICMs over $K_2$ was achieved in \cite{CL}.
In particular, the rank of a nontrivial FICM over $K_2$ is either $(1+1)$ or $(2+2)$.
Obviously, any conformal module over $K_2\subset\KB$ can be trivially extended to a conformal module over $\KB$.
Our main result indicates that a FICM over $\KB$ admits a nontrivial
extension of a finite conformal module $M$ over $K_2$ if and only if $p=-1$ and $M$ has rank $(2+2)$ (see Table~1).
As a byproduct of our main result, we also obtain the classification of FICMs over the finite Lie conformal superalgebra ${\frak k}(n)$ (see Table~2).

\begin{table}[h]
\centering\small
\subtable[Nontrivial FICMs over $\KB$]{
\begin{tabular}{c|c|c}
\hline
 $\KB$  & FICMs \footnotesize{(up to parity change)} & Reference \\
\hline
$p\ne -1$  & $V^{(1)}_{\D,\a}$, $V^{(2)}_{\D,\a}$, $V_{\D,\L,\a}$ & Theorem \ref{thm-classification} \\
\hline
$p=-1$ & $V^{(1)}_{\D,\a}$, $V^{(2)}_{\D,\a}$, $V_{\D,\L,\a,\b}$ & Theorem \ref{thm-classification}\\
\hline
\end{tabular}}
\\[10pt]
\subtable[Nontrivial FICMs over ${\frak k}(n)$]{
\begin{tabular}{c|c|c}
\hline
${\frak k}(n)$  & FICMs \footnotesize{(up to parity change)} & Reference \\
\hline
$n> 1$  & $V^{(1)}_{\D,\a}$, $V^{(2)}_{\D,\a}$, $V_{\D,\L,\a}$ & Corollary~\ref{classification-kn} \\
\hline
$n=1$ & $V^{(1)}_{\D,\a}$, $V^{(2)}_{\D,\a}$, $V_{\D,\L,\a,\b}$ & Corollary~\ref{classification-kn}\\
\hline
\end{tabular}}
\end{table}

To our surprising, any rank $(1+1)$ $K_2$-module $M$
can not be nontrivially extended to a $\KB$-module, even in case that $p=-1$ and $M$ can be degenerated from a rank $(2+2)$ $K_2$-module (see Remark~\ref{remark-1+1-triviality}).
This is essentially different from the extensions from $K_1$-modules to $\SB$-modules \cite{X}.
As one of key steps to prove our main results, we classify all the free conformal modules
of ranks $(1+1)$ and $(2+2)$ over $\KB$ and ${\frak k}(n)$ without the irreducibility assumption,
and completely determine their composition factors (see Tables~3 and 4 in Proposition~\ref{prop-simplicity-2+2}).
As a corollary, the composition factors of $K_2$-modules of small rank are also obtained (see Remark~\ref{remark-K2-CF}).

The outline of the paper is as follows.
In Section~2, after recalling some basic definitions and notations, we construct the main object $\KB$,
and present certain quotient algebras $\KB_{[n]}$ (cf.~\eqref{equ-quotiont}) of $\KB$.
In particular, we introduce the $N=2$ Heisenberg and Schr$\ddot{\rm o}$dinger conformal superalgebras.
In Section~3, for certain subquotient algebra ${\frak g}$ (cf.~\eqref{subquotient}) of the annihilation superalgebra of $\KB$,
by introducing the {\it row ideals} and {\it column ideals} (cf.~\eqref{r-ideal} and \eqref{c-ideal}),
we conceptually determine the dimension of any irreducible module over ${\frak g}$.
In Section~4, we classify all the free conformal $\KB$-modules
of ranks $(1+1)$ and $(2+2)$ by analytical techniques.
Also, for reducible ones, we completely determine their composition factors.
In the last section, we complete the classification of FICMs over $\KB$ by showing that they must be free of rank $(1+1)$ or $(2+2)$.
As an application, we also obtain the classification of FICMs over ${\frak k}(n)$ by the feature \eqref{equ-sn} of $\KBB(-n)$.

\section{\large{Preliminaries}}

Throughout this paper, the notation $|a|\in\Z/2\Z$ denotes the parity of an element $a$ in a super vector space,
and $a$ is always assumed to be homogeneous if $|a|$ appears in an expression.
The angle bracket $\langle\ ,\rangle$ denotes ``the Lie conformal superalgebra generated over $\C[\partial]$ by''.
The symbol $\delta$ denotes the Kronecker delta.

\subsection{Basic definitions}
Let us first recall some basic definitions, see \cite{FK,K1,X} for more details.

\begin{defi}\label{def-conforma-superalgebra}\rm
A {\it Lie conformal superalgebra} $R=R_{\bar{0}}\oplus R_{\bar{1}}$ is a $\Z/2\Z$-graded $\C[\partial]$-module endowed with a
$\C$-linear map $R\otimes R\rightarrow \C[\lambda]\otimes R$, $a\otimes b\rightarrow [a\,{}_\lambda \,b]$
called $\lambda$-bracket, and satisfying the following axioms ($a,\,b,\,c\in R$):
\begin{equation*}
\aligned
\mbox{(conformal sesquilinearity)}&~~~~[\partial a\,{}_\lambda \,b]=-\lambda[a\,{}_\lambda\, b],\ \ \ \
[a\,{}_\lambda \,\partial b]=(\partial+\lambda)[a\,{}_\lambda\, b],\\
\mbox{(skew-symmetry)}&~~~~[a\, {}_\lambda\, b]=-(-1)^{|a||b|}[b\,{}_{-\lambda-\partial}\,a],\\
\mbox{(Jacobi identity)}&~~~~[a\,{}_\lambda\,[b\,{}_\mu\, c]]=[[a\,{}_\lambda\, b]\,{}_{\lambda+\mu}\, c]+(-1)^{|a||b|}[b\,{}_\mu\,[a\,{}_\lambda \,c]].
\endaligned
\end{equation*}
\end{defi}

Let $R$ be a Lie conformal superalgebra. We call $R$ {\it finite} if it is finitely generated over $\C[\partial]$;
{\it $\Z$-graded} if $R=\oplus_{i\in \Z}R_i$, where $R_i$ is a $\C[\partial]$-submodule
and $[R_i\,{}_\lambda\,R_j]\subset R_{i+j}[\lambda]$ for $i,j\in \Z$.

\begin{defi}\label{def-conformal-module}\rm
A {\it conformal module} $M=M_{\bar{0}}\oplus M_{\bar{1}}$ over a Lie conformal superalgebra $R$ is a $\Z/2\Z$-graded $\C[\partial]$-module endowed with a $\C$-linear map
$R\otimes M\rightarrow \C[\lambda]\otimes M$, $a\otimes v\rightarrow a\,{}_\lambda \,b$ called $\lambda$-action, such that ($a,\,b\in R$, $v\in M$)
\begin{equation*}
\aligned
&(\partial a)\,{}_\lambda\, v=-\lambda a\,{}_\lambda\, v,\ \ \ \ \ a{}\,{}_\lambda\, (\partial v)=(\partial+\lambda)a\,{}_\lambda\, v,\\
&[a\,{}_\lambda\, b]\,{}_{\lambda+\mu}\, v = a\,{}_\lambda\, (b{}\,_\mu\, v)-(-1)^{|a||b|}b\,{}_\mu\,(a\,{}_\lambda\, v).
\endaligned
\end{equation*}
\end{defi}

Let $M=M_{\bar{0}}\oplus M_{\bar{1}}$ be a conformal $R$-module.
We call $M$ {\it finite} if it is finitely generated over $\C[\partial]$.
As $\C[\partial]$-modules, if $M_{\bar{0}}$ has rank $m$ and $M_{\bar{1}}$ has rank $n$, we say that $M$ has {\it rank $(m+n)$}, denoted by ${\rm{rank}}(M)=(m+n)$.
In case $R$ is $\Z$-graded, we call $M$ {\it $\Z$-graded} if $M=\oplus_{i\in \Z}M_i$, where $M_i$ is a $\C[\partial]$-submodule and
$R_i\,{}_\lambda\, M_j\subset M_{i+j}[\lambda]$ for $i,j\in \Z$. Furthermore, if each $M_i$ is freely generated by one element over $\C[\partial]$,
we call $M$ a {\it $\Z$-graded free intermediate series module}.

Clearly, for any fixed $\a\in\C$, the $\C[\partial]$-module $\C c_\a$ with $\partial c_\a=\a c_\a$, $a\,{}_\lambda c_\a=0$ for $a\in R$, is a conformal $R$-module,
which will be referred to as the {\it even} (resp., {\it odd}) {\it one-dimensional trivial module} if $|c_\a|=\bar{0}$ (resp., $|c_\a|=\bar{1}$).

\begin{defi}\label{def-annihilation-algebra}\rm
The {\it  annihilation superalgebra} ${\cal A}(R)$ of a Lie conformal superalgebra $R$ is a Lie superalgebra with
$\C$-basis $\{a_{n}\,|\,a\in R,\,n\in\Z_+\}$ and relations
\begin{equation}\label{equ-new}
\aligned
&(\lambda a)_{n}=\lambda a_{n},\quad (a+b)_{n}=a_{n}+b_{n},\quad (\partial a)_{n}=-n a_{n-1},\\
&[a_{m}, b_{n}]=\sum_{k\in\Z_+}{\binom m k}(a_{(k)}b)_{m+n-k},
\endaligned
\end{equation}
where $a_{(k)}b$ is called the $k$-product, given by the following inversion formula:
$$
[a\,{}_\lambda \, b]=\sum_{k\in\Z_+}\lambda^{(k)}a_{(k)}b\quad\mbox{with}\quad \lambda^{(k)}=\frac{\lambda^k}{k!}.
$$
\end{defi}

Here, the parity $|a_{n}|$ of $a_{n}\in{\cal A}(R)$ is the same as $|a|$ for any $a\in R$ and $n\in\Z_+$.
Note that ${\cal A}(R)$ admits a derivation $T$ given by $T(a_{n})=-n a_{n-1}$ for any $a_{n}\in {\cal A}(R)$.
The {\it extended annihilation superalgebra} ${\cal A}(R)^e$ of a Lie conformal superalgebra $R$ is defined by
${\cal A}(R)^e=\C T\ltimes{\cal A}(R)$ with $[T, a_{n}]=-n a_{n-1}$.
The representation theory of $R$ is controlled by the representation theory of ${\cal A}(R)^e$ in the following sense:

\begin{prop}\label{prop-observation}
A conformal module $M$ over a Lie conformal superalgebra $R$ is the same as
a module over the Lie superalgebra ${\cal A}(R)^e$ satisfying $a_{n}v=0$ for $a\in R$, $v\in M$, $n\gg 0$.
\end{prop}

\subsection{Construction of $\KB$}

Recall that \cite{FK} the $N=2$ conformal superalgebra $K_2$ is a Cartan $K$ type Lie conformal superalgebra,
which has $\C[\partial]$-basis $\{L,\,J,\,G^{\pm}\}$ with $|L|=|J|=\bar{0}$ and $|G^{\pm}|=\bar{1}$, satisfying
\begin{eqnarray}
\nonumber[L\, {}_\lambda \, L] &\!\!\!=\!\!\!& (\partial+2\lambda) L, \\[-2pt]
\nonumber[L\, {}_\lambda \, J] &\!\!\!=\!\!\!& (\partial+\lambda) J, \\[-2pt]
\nonumber[L\, {}_\lambda \, G^{\pm}] &\!\!\!=\!\!\!& (\partial+\frac{3}{2}\lambda) G^{\pm}, \\[-2pt]
\nonumber[J\, {}_\lambda \, G^{\pm}] &\!\!\!=\!\!\!& \pm G^{\pm}, \\[-2pt]
\nonumber[G^{+}\, {}_\lambda \, G^{-}] &\!\!\!=\!\!\!& (\partial+2\lambda) J+2 L.
\end{eqnarray}
The Lie conformal algebra $\CB$ defined by \eqref{brackets-LL} was introduced in \cite{SXY}.
Motivated by an important class of $\Z$-graded free intermediate series $\CB$-module structure \cite{X}
and the conformal structure of $K_2$, let us consider a $\Z/2\Z$-graded $\C[\partial]$-module
$$
R:=R\left(\left\{\phi^{(k)}_{i,j}|\,k=1,2,3,4\right\}\right)=R_{\bar{0}}\oplus R_{\bar{1}},
$$
where the even part $R_{\bar{0}}$ has $\C[\partial]$-basis
$\{L_i,\,J_i\,|\,i\in\Z_+\}$ and the odd part $R_{\bar{1}}$ has $\C[\partial]$-basis $\{G^{\pm}_i\,|\,i\in\Z_+\}$,
satisfying \eqref{brackets-LL}--\eqref{brackets-LG} and
\begin{equation*}
\aligned
&[J_i\, {}_\lambda \, G^+_j] = \phi^{(1)}_{i,j}(\partial,\lambda) G^+_{i+j}, \quad\quad [J_i\, {}_\lambda \, G^-_j] = \phi^{(2)}_{i,j}(\partial,\lambda) G^-_{i+j},\\
&[G^+_i\, {}_\lambda \, G^-_j] = \phi^{(3)}_{i,j}(\partial,\lambda) J_{i+j} + \phi^{(4)}_{i,j}(\partial,\lambda) L_{i+j},\quad [X_i\, {}_\lambda \, X_j]=0\mbox{ for }X=J,\,G^{\pm}.
\endaligned
\end{equation*}
Here, $\phi^{(k)}_{i,j}(\partial,\lambda)\in\C[\partial,\lambda]$, $k=1,2,3,4$, referred to as the {\it structure polynomials},
satisfy $\phi^{(k)}_{i,j}(\partial,\lambda)\ne0$ for some $i,j\in\Z_+$.

\begin{lemm}\label{construction}
The $\Z/2\Z$-graded $\C[\partial]$-module $R$ becomes a Lie conformal superalgebra if and only if
$\phi^{(1)}_{i,j}(\partial,\lambda)=-\phi^{(2)}_{i,j}(\partial,\lambda)=a$,
$\phi^{(3)}_{i,j}(\partial,\lambda)=\frac{c}{a}((i+\frac{p}{2})\partial+(i+j+p)\lambda)$, $\phi^{(4)}_{i,j}(\partial,\lambda)=c$
for all $i,j\in\Z_+$, where $a$ and $c$ are constants.
Up to isomorphism, we may assume that $a=1,\,c=2$,
and the resulting algebra is exactly $\KB$.
\end{lemm}

\ni{\it Proof.}\ \
The sufficiency is straightforward, so let us prove the necessity.
Let $R$ be a Lie conformal superalgebra. We need to determine the structure polynomials $\phi^{(k)}_{i,j}(\partial,\lambda)$, $k=1,2,3,4$.
Roughly speaking, our strategy is to apply the Jacobi identity for certain triples, and then
to determine $\phi^{(k)}_{i,j}(\partial,\lambda)$ by using the substitutability of variables, the arbitrariness of subscripts and some analytical techniques.

Using the Jacobi identity for triple $(L_0,G_i^+,G_j^-)$ and matching the coefficients of $J_{i+j}$ and $L_{i+j}$, we obtain
\begin{eqnarray*}
(p\partial+(i+j+p)\lambda)\phi_{i,j}^{(3)}(\partial+\lambda,\mu)&\!\!\!=\!\!\!&
((i+\frac{1}{2}p)\lambda-p\mu)\phi_{i,j}^{(3)}(\partial,\lambda+\mu)\\
&\!\!\!\!\!\!&
+(p(\partial+\mu)+(j+\frac{3}{2}p)\lambda)\phi_{i,j}^{(3)}(\partial,\mu),\\
(p\partial+(i+j+2p)\lambda)\phi_{i,j}^{(4)}(\partial+\lambda,\mu)&\!\!\!=\!\!\!&
((i+\frac{1}{2}p)\lambda-p\mu)\phi_{i,j}^{(4)}(\partial,\lambda+\mu)\\
&\!\!\!\!\!\!&
+(p(\partial+\mu)+(j+\frac{3}{2}p)\lambda)\phi_{i,j}^{(4)}(\partial,\mu),
\end{eqnarray*}
where $i,j\in\Z_+$.
Taking $\mu=0$, we have
\begin{eqnarray*}
(i+\frac{1}{2}p)\phi_{i,j}^{(3)}(\partial,\lambda)=
p\partial\frac{\phi_{i,j}^{(3)}(\partial+\lambda,0)-\phi_{i,j}^{(3)}(\partial,0)}{\lambda}+(i+j+p)\phi_{i,j}^{(3)}(\partial+\lambda,0)-(j+\frac{3}{2}p)\phi_{i,j}^{(3)}(\partial,0),\\
(i+\frac{1}{2}p)\phi_{i,j}^{(4)}(\partial,\lambda)=
p\partial\frac{\phi_{i,j}^{(4)}(\partial+\lambda,0)-\phi_{i,j}^{(4)}(\partial,0)}{\lambda}+(i+j+2p)\phi_{i,j}^{(4)}(\partial+\lambda,0)-(j+\frac{3}{2}p)\phi_{i,j}^{(4)}(\partial,0).
\end{eqnarray*}
Taking $\lambda\rightarrow0$, we obtain
\begin{equation}\label{poly34}
\partial \frac{d}{d\partial}(\phi_{i,j}^{(3)}(\partial,0))=\phi_{i,j}^{(3)}(\partial,0),\quad
\partial \frac{d}{d\partial}(\phi_{i,j}^{(4)}(\partial,0))=0.
\end{equation}
The second formula in \eqref{poly34} implies that $\phi_{i,j}^{(4)}(\partial,0)=c_{i,j}\in\C$.
Processing with the same procedure for triple $(L_i, G_0^+, G_j^-)$ and matching the coefficients of $L_{i+j}$, we obtain
\begin{equation}\label{poly34-4}
\lambda\left(\frac{1}{2}p\phi_{i,j}^{(4)}(\partial,\lambda)+(i+j+\frac{3}{2}p)c_{0,i+j}-(i+j+2p)c_{0,j}\right)=
(i+p)\partial(c_{0,j}-c_{0,i+j}).
\end{equation}
By the arbitrariness of $i$ and $j$, we have $\lambda|(c_{0,j}-c_{0,i+j})$, and so $c_{0,j}=c$ for all $j\in\Z_+$, where $c$ is a constant.
Substituting this back to \eqref{poly34-4}, we see that $\phi_{i,j}^{(4)}(\partial,\lambda)=c$ for all $i,j\in\Z_+$.

To save space, we proceed with an outline of the characterization of the other structure polynomials; the details are omitted.
For $\phi^{(k)}_{i,j}(\partial,\lambda)$, $k=1,2$, by considering triple $(L_m,J_i,G_j^\pm)$, one can
first show that $\phi_{i,j}^{(k)}(\partial,\lambda)=a^{(k)}$ for all $i,j\in\Z_+$, where $a^{(k)}$ is a constant.
Then, by triple $(G_0^+,J_i,G_j^-)$ one can derive that $a^{(1)}=-a^{(2)}$ (simply denote $a^{(1)}$ by $a$) and obtain a preliminary
form of $\phi^{(3)}_{i,j}(\partial,\lambda)$. Using this preliminary form, the first formula in \eqref{poly34}, and triple $(G_0^+,G_0^+,G_j^-)$,
one can prove that $\phi^{(3)}_{i,j}(\partial,\lambda)$ must have the required form.

For the last statement, we only need to note that the map $L_i\mapsto L_i$, $J_i\mapsto\frac{1}{a}J_i$, $G_i^\pm\mapsto\sqrt{\frac{2}{c}}G_i^\pm$
defines an automorphism of $R$. This completes the proof.
\QED
\vskip5pt

\subsection{Quotient algebras of $\KB$}

One of attractive features of $\KB$, as we mentioned in the Introduction, is that one can obtain new interesting
finite Lie conformal superalgebras from $\KB$.
Let $n$ be a positive integer. Consider a subalgebra $\KB_{\langle n\rangle}$ of $\KB$ defined by
\begin{equation*}\label{equ-ideal}
\KB_{\langle n\rangle}= \langle L_i,\,J_i,\,G_i^\pm\,|\,i\ge n \rangle.
\end{equation*}
Note that $\KB_{\langle n\rangle}$ is in fact a Lie conformal superalgebra ideal of $\KB$.
For $n\in\Z_+$, define a quotient algebra $\KB_{[n]}$ of $\KB$ by
\begin{equation}\label{equ-quotiont}
\KB_{[n]}=\KB/\KB_{\langle n+1\rangle}.
\end{equation}
Clearly, $\KB_{[0]}\cong K_2$.
In addition, the special cases $p=-n$ with $n\in \Z_{\ge 1}$ supply a series of  new finite non-simple Lie conformal superalgebras:
\begin{equation}\label{equ-quotiont-special}
{\frak k}(n)=\KBB(-n)_{[n]}=\KBB(-n)/\KBB(-n)_{\langle n+1 \rangle}.
\end{equation}
Let us say more about the conformal structures of ${\frak k}(1)$ and ${\frak k}(2)$.

\begin{exam}\rm
By definition \eqref{equ-quotiont-special}, one can check that the subalgebra $\langle \bar{L}_0, \bar{J}_0, \bar{G}_0^\pm \rangle$ of ${\frak k}(1)$
is isomorphic to the $N=2$ conformal superalgebra $K_2$.
While the subalgebra $\langle \bar{L}_i, \bar{G}_i^+ +\bar{G}_i^-\,|\,i=0,1\rangle$ of ${\frak k}(1)$
is isomorphic to the Heisenberg-Neveu-Schwarz conformal algebra \cite{X}.
Following our previous name rule, we refer to ${\frak k}(1)$ as {\it $N=2$ Heisenberg conformal superalgebra}.
\end{exam}

\begin{exam}\rm
Similarly, the subalgebra $\langle \bar{L}_0, \bar{J}_0, \bar{G}_0^\pm \rangle$ of ${\frak k}(2)$
is also isomorphic to the $N=2$ conformal superalgebra $K_2$.
While the subalgebra $\langle \bar{L}_i, \bar{G}_i^+ +\bar{G}_i^-\,|\,i=0,1,2\rangle$ of ${\frak k}(2)$
is isomorphic to the Schr$\ddot{\rm o}$dinger-Neveu-Schwarz conformal algebra \cite{X}.
Following our previous name rule, we refer to ${\frak k}(2)$ as {\it $N=2$ Schr$\ddot{\rm o}$dinger conformal superalgebra}.
\end{exam}

\section{\large{Annihilation superalgebra and related representations}}

In this section, we construct a subquotient algebra ${\frak g}$ (cf.~\eqref{subquotient}) of the annihilation superalgebra ${\cal A}(\KB)$ of $\KB$.
Then, by introducing the so-called row ideals and column ideals of ${\frak g}$, we determine the dimension of any nontrivial finite-dimensional irreducible module over ${\frak g}$.

\subsection{Annihilation superalgebra ${\cal A}(\KB)$}

Let us first give the Lie superalgebra structure on the annihilation superalgebra ${\cal A}(\KB)$.

\begin{lemm}\label{annihilation-algebra}
We have
$$
{\cal A}(\KB)\cong\sp_{\C} \{L_{i,m},\,J_{i,n},\,G_{i,t}^\pm\,|\,i,n\in\Z_{+},\,m\in\Z_{\ge-1},\,t\in\frac{1}{2}+\Z_{\ge-1}\}.
$$
with super-commutation relations
\begin{eqnarray*}
\nonumber [L_{i,m},L_{j,n}] &\!\!\!=\!\!\!& ((j+p)(m+1)-(i+p)(n+1)) L_{i+j,m+n},\\
\nonumber [L_{i,m},J_{j,n}] &\!\!\!=\!\!\!& (j(m+1)-(i+p)n) J_{i+j,m+n},\\
\nonumber [L_{i,m},G_{j,n}^\pm] &\!\!\!=\!\!\!& ((j+\frac{p}{2})(m+1)-(i+p)(n+\frac{1}{2})) G_{i+j,m+n}^\pm,\\
\nonumber [J_{i,m},G_{j,n}^\pm] &\!\!\!=\!\!\!& \pm G_{i+j,m+n}^\pm,\\
\nonumber [G_{i,m}^+,G_{j,n}^-] &\!\!\!=\!\!\!& ((2j+p)(m+\frac{1}{2})-(2i+p)(n+\frac{1}{2}))J_{i+j,m+n}+2 L_{i+j,m+n}.
\end{eqnarray*}
\end{lemm}

\ni{\it Proof.}\ \
By the inversion formula in Definition~\ref{def-annihilation-algebra}, one can first transfer the
$\lambda$-brackets of $\KB$ to $k$-products:
\begin{equation*}
\aligned
&{L_i}_{(k)}{L_j}=
\left\{\begin{array}{ll}
(i+p)\partial L_{i+j}, &\mbox{if \ } k=0,\\[4pt]
(i+j+2p)L_{i+j}, &\mbox{if \ } k=1,\\[4pt]
0, &\mbox{if \ } k\ge 2,\end{array}\right.
\quad
{L_i}_{(k)}{J_j}=
\left\{\begin{array}{ll}
(i+p)\partial J_{i+j}, &\mbox{if \ } k=0,\\[4pt]
(i+j+p)J_{i+j}, &\mbox{if \ } k=1,\\[4pt]
0, &\mbox{if \ } k\ge 2,\end{array}\right.\\
&{L_i}_{(k)}{G_j^\pm}=
\left\{\begin{array}{ll}
(i+p)\partial G_{i+j}^\pm, &\mbox{if \ } k=0,\\[4pt]
(i+j+\frac{3}{2}p)G_{i+j}^\pm, &\mbox{if \ } k=1,\\[4pt]
0, &\mbox{if \ } k\ge 2,\end{array}\right.
\quad
{J_i}_{(k)}{G_j^\pm}=
\left\{\begin{array}{ll}
\pm G_{i+j}^\pm, &\mbox{if \ } k=0,\\[4pt]
0, &\mbox{if \ } k\ge 1,\end{array}\right.\\
&{G_i^+}_{(k)}{G_j^-}=
\left\{\begin{array}{ll}
(2i+p)\partial J_{i+j} + 2L_{i+j}, &\mbox{if \ } k=0,\\[4pt]
2(i+j+p)J_{i+j}, &\mbox{if \ } k=1,\\[4pt]
0, &\mbox{if \ } k\ge 2.\end{array}\right.
\endaligned
\end{equation*}
Then, by \eqref{equ-new} and making the shift $L_{i,m}=(L_i)_{m+1}$, $J_{i,n}=(J_i)_{n}$, $G_{i,t}^\pm=(G_i^\pm)_{t+\frac{1}{2}}$
for $i,n\in\Z$, $m\in\Z_{\ge-1}$, $t\in\frac{1}{2}+\Z_{\ge-1}$,
one can obtain the super-commutation relations of ${\cal A}(\KB)$.
\QED

\subsection{Subquotient algebras of ${\cal A}(\KB)$}

Next, we construct a class of subquotient algebras of ${\cal A}(\KB)$.
Let $k,N\in\Z_+$ be two non-negative integers.
Consider a subalgebra ${\cal A}(\KB)_+$ of ${\cal A}(\KB)$:
$$
{\cal A}(\KB)_+=\sp_{\C}\{L_{i,m},\,J_{i,n},\,G_{i,t}^\pm\in{\cal A}(\KB)\,|\,i,m\in\Z_{+},\,n\in\Z_{\ge1},\,t\in\frac{1}{2}+\Z_{+}\}.
$$
This Lie superalgebra contains an ideal ${\cal I}(k,N)$:
$$
{\cal I}(k,N)=\sp_{\C}\{L_{i,m},\,J_{i,n},\,G_{i,t}^\pm\in{\cal A}(\KB)_+\,|\, i>k, \, m> N,\,n>N+1,\,t> N+\frac{1}{2}\}.
$$
Define ${\frak g}(k,N)$ by
\begin{equation}\label{subquotient}
{\frak g}(k,N):={\cal A}(\KB)_+/{\cal I}(k,N).
\end{equation}
We shall determine the dimension of any finite-dimensional irreducible module over ${\frak g}(k,N)$.
To this end, we introduce two types of ideals of ${\frak g}(k,N)$:
\begin{itemize}
\item {\it row ideals} $\mathfrak{r}(k,N)$ defined by
\begin{equation}\label{r-ideal}
\mathfrak{r}(k,N) = \sp_{\C}\{\bar{L}_{k,n},\,\bar{J}_{k,n+1},\,\bar{G}_{k,n+\frac{1}{2}}^\pm\in{\frak g}(k,N)\,|\, n\le N\};
\end{equation}
\item {\it column ideals} $\mathfrak{c}(k,N)$ defined by
\begin{equation}\label{c-ideal}
\mathfrak{c}(k,N) = \sp_{\C}\{\bar{L}_{i,N},\,\bar{J}_{i,N+1},\,\bar{G}_{i,N+\frac{1}{2}}^\pm\in{\frak g}(k,N)\,|\, i\le k\}.
\end{equation}
\end{itemize}

\subsection{Irreducible representations of ${\frak g}(k,N)$}

Let $V=V_{\bar{0}}\oplus V_{\bar{1}}$ be a nontrivial finite-dimensional irreducible module over ${\frak g}(k,N)$.
If $V_{\bar{0}}$ has dimension $m$ and $V_{\bar{1}}$ has dimension $n$, we say that $V$ has dimension $(m|n)$.
The main result in this section is the following.

\begin{theo}\label{thm-1}
The dimension of $V$ is $(1|1)$ or $(2|2)$.
\end{theo}

Let us first introduce some auxiliary sets: 
\begin{eqnarray}
\nonumber \Omega &\!\!\!=\!\!\!& \{(j,n)\,|\, \bar{L}_{j,n},\,\bar{J}_{j,n+1},\,\bar{G}_{j,n+\frac{1}{2}}^\pm\in{\frak g}(k,N)\}\backslash\{(0,0)\},\\
\label{aux-L} \Omega(L) &\!\!\!=\!\!\!& \{(j,n)\in\Omega \,|\, j-pn=0\},\\
\label{aux-J} \Omega(J) &\!\!\!=\!\!\!& \{(j,n)\in\Omega \,|\, j-p(n+1)=0\},\\
\label{aux-G} \Omega(G) &\!\!\!=\!\!\!& \{(j,n)\in\Omega \,|\, j-p(n+\frac{1}{2})=0\}.
\end{eqnarray}

The following facts can be easily checked.

\begin{prop}\label{prop-fact}
Suppose $k,\,N\ge1$. We have
\baselineskip1pt\lineskip7pt\parskip-1pt
\begin{itemize}\parskip-1pt
  \item[{\rm(1)}] $\Omega(L)\ne\emptyset \Longleftrightarrow \Omega(J)\ne\emptyset$, and if $\Omega(L)\ne\emptyset$, then $p>0$;
  \item[{\rm(2)}] $\Omega(G)\ne\emptyset \Longleftrightarrow p\in 2\Z_{\ge1}$, and if $\Omega(G)\ne\emptyset$, then $\Omega(L)\ne\emptyset$.
\end{itemize}
\end{prop}

\begin{lemm}\label{lemma-1}
If $\Omega(L)=\Omega(J)=\Omega(G)=\emptyset$, then $\dim V=(1|1)$.
\end{lemm}

\ni{\it Proof.}\ \
Consider the following decomposition of ${\frak g}(k,N)$:
$$
{\frak g}(k,N)=\C \bar{L}_{0,0}+\check{{\frak g}}(k,N), \mbox{\ \ where\ \ }\check{{\frak g}}(k,N)={\frak g}(k,N)\backslash\C \bar{L}_{0,0}.
$$
Clearly, $\check{{\frak g}}(k,N)$ is a nilpotent ideal of ${\frak g}(k,N)$.
Since $\Omega(L)=\Omega(J)=\Omega(G)=\emptyset$, one can see that $\check{{\frak g}}(k,N)$
is a completely reducible $\C \bar{L}_{0,0}$-module with no trivial summand.
By \cite[Lemma~1]{CK}, $\check{{\frak g}}(k,N)$ acts trivially on $V$. Hence, $V$ is simply a finite-dimensional $\C \bar{L}_{0,0}$-module,
and so $\dim V=(1|1)$.
\QED

\begin{lemm}\label{lemma-2}
Suppose $k,\,N\ge1$ and $\Omega(L)\ne \emptyset$, $\Omega(G)=\emptyset$. Let
$$
\begin{aligned}
&  j_L={\rm{max}}\{j\,|\,(j,n)\in\Omega(L)\}, && n_L={\rm{max}}\{n\,|\,(j,n)\in\Omega(L)\}, \\[-2pt]
&  j_J={\rm{max}}\{j\,|\,(j,n)\in\Omega(J)\}, && n_J={\rm{max}}\{n\,|\,(j,n)\in\Omega(J)\}.
\end{aligned}
$$
\baselineskip1pt\lineskip7pt\parskip-1pt
\begin{itemize}\parskip-1pt
  \item[{\rm(1)}] If $j_J<k$, then the row ideal $\mathfrak{r}(k,N)$ acts trivially on $V$;
  \item[{\rm(2)}] If $j_J=k$, $j_L<k$, then the row ideal $\mathfrak{r}(k,N)$ acts trivially on $V$;
  \item[{\rm(3)}] If $j_J=j_L=k$, $n_L<N$, then the column ideal $\mathfrak{c}(k,N)$ acts trivially on $V$;
  \item[{\rm(4)}] If $j_J=j_L=k$, $n_L=N$, then the column ideal $\mathfrak{c}(k,N)$ acts trivially on $V$.
\end{itemize}
\end{lemm}

\ni{\it Proof.}\ \
Note first that $j_L\le j_J$ and (1)--(4) cover all possible cases.

(1) Assume that $\mathfrak{r}(k,N)$ acts non-trivially on $V$. We have
$V=\mathfrak{r}(k,N)V$ by the irreducibility of $V$.
Consider the action of $\bar{L}_{0,0}$ on $\mathfrak{r}(k,N)$,
\begin{eqnarray}
\label{equ-L0L} [\bar{L}_{0,0}, \bar{L}_{k,n}] &\!\!=\!\!& b_1\bar{L}_{k,n}, \mbox{\ where\ } b_1=k-pn\ne0 \ (\mbox{since\ } k>j_J\ge j_L), \\
\label{equ-L0J} [\bar{L}_{0,0}, \bar{J}_{k,n+1}] &\!\!=\!\!& b_2\bar{J}_{k,n+1}, \mbox{\ where\ } b_2=k-p(n+1)\ne0 \ (\mbox{since\ } k>j_J), \\
\label{equ-L0G} [\bar{L}_{0,0}, \bar{G}_{k,n+\frac{1}{2}}^\pm] &\!\!=\!\!& b_3\bar{G}_{k,n+\frac{1}{2}}^\pm, \mbox{\ where\ } b_3=k-p(n+\frac{1}{2})\ne0 \ (\mbox{since\ } \Omega(G)=\emptyset).
\end{eqnarray}
Hence, $\bar{L}_{k,n}, \bar{J}_{k,n+1}, \bar{G}_{k,n+\frac{1}{2}}^\pm\in\mathfrak{r}(k,N)$ act nilpotently on $V$.
Since $\mathfrak{r}(k,N)$ is abelian, $\mathfrak{r}(k,N)$ acts nilpotently on $V$, which contradicts to $V=\mathfrak{r}(k,N)V$.

(2) In this case, we have $n_J=N$.
Assume that $\mathfrak{r}(k,N)$ acts non-trivially on $V$. We have
$V=\mathfrak{r}(k,N)V$ by the irreducibility of $V$.
Consider the decomposition of $\mathfrak{r}(k,N)$:
$$
{\frak r}(k,N)=\C \bar{J}_{k,N+1}+\check{{\frak r}}(k,N), \mbox{\ \ where\ \ }\check{{\frak r}}(k,N)={\frak r}(k,N)\backslash\C \bar{J}_{k,N+1}.
$$
Note that $\bar{J}_{k,N+1}$ is an even central element of $\mathfrak{g}(k,N)$. We may assume that the action of $\bar{J}_{k,N+1}$ is a scalar $c$.
In addition, we have
\begin{equation}\label{L-case2}
[\bar{L}_{k,0}, \bar{J}_{0,N+1}] = -(N+1)(k+p)\bar{J}_{k,N+1}.
\end{equation}
Consider the actions of both sides of \eqref{L-case2} on $V$, and compare the traces of the matrices with respect to a basis of $V$.
The right hand side is $-c(N+1)(k+p)\dim V$, while the left hand side is zero, since the corresponding matrix has form $AB-BA$.
Hence, $c=0$ (note that $k+p>0$ by Proposition~\ref{prop-fact}(1)), and so $V=\check{{\frak r}}(k,N)V$.
Considering the action of $\bar{L}_{0,0}$ on $\check{{\frak r}}(k,N)$, we still have \eqref{equ-L0L}--\eqref{equ-L0G}.
The only difference is the reason for $b_2\ne0$, which should be replaced by $n<N=n_J$.
Hence, all elements in $\check{{\frak r}}(k,N)$ act nilpotently on $V$. Since $\check{{\frak r}}(k,N)$ is abelian, $\check{{\frak r}}(k,N)$
acts nilpotently on $V$, which contradicts to $V=\check{{\frak r}}(k,N)V$.

(3) In this case, we have $n_J=n_L-1$.
Assume that $\mathfrak{c}(k,N)$ acts non-trivially on $V$. We have
$V=\mathfrak{c}(k,N)V$ by the irreducibility of $V$.
Consider the action of $\bar{L}_{0,0}$ on $\mathfrak{c}(k,N)$,
\begin{eqnarray}
\label{equ-L0L-c} [\bar{L}_{0,0}, \bar{L}_{i,N}] &\!\!=\!\!& b_4\bar{L}_{i,N}, \mbox{\ where\ } b_4=i-pN\ne0 \ (\mbox{since\ } N>n_L), \\
\label{equ-L0J-c} [\bar{L}_{0,0}, \bar{J}_{i,N+1}] &\!\!=\!\!& b_5\bar{J}_{i,N+1}, \mbox{\ where\ } b_5=i-p(N+1)\ne0 \ (\mbox{since\ } N>n_L>n_J), \\
\label{equ-L0G-c} [\bar{L}_{0,0}, \bar{G}_{i,N+\frac{1}{2}}^\pm] &\!\!=\!\!& b_6\bar{G}_{i,N+\frac{1}{2}}^\pm, \mbox{\ where\ } b_6=i-p(N+\frac{1}{2})\ne0 \ (\mbox{since\ } \Omega(G)=\emptyset).
\end{eqnarray}
Hence, $\bar{L}_{i,N}, \bar{J}_{i,N+1}, \bar{G}_{i,N+\frac{1}{2}}^\pm\in\mathfrak{c}(k,N)$ act nilpotently on $V$. Since $\mathfrak{c}(k,N)$ is abelian, $\mathfrak{c}(k,N)$
acts nilpotently on $V$, which contradicts to $V=\mathfrak{c}(k,N)V$.

(4) In this case, we have $n_J=N-1$.
Assume that $\mathfrak{c}(k,N)$ acts non-trivially on $V$.
Consider the decomposition of ${\frak c}(k,N)$:
\begin{eqnarray*}
{\frak c}(k,N) &\!\!\!=\!\!\!& \mathfrak{c}_L(k,N)+\mathfrak{c}_{JG}(k,N), \quad \mbox{where} \\
\mathfrak{c}_L(k,N) &\!\!\!=\!\!\!& \sp_{\C}\{\bar{L}_{i,N}\in{\frak g}(k,N)\,|\, i\le k\}, \\
\mathfrak{c}_{JG}(k,N) &\!\!\!=\!\!\!& \sp_{\C}\{\bar{J}_{i,N+1},\,\bar{G}_{i,N+\frac{1}{2}}^\pm\in{\frak g}(k,N)\,|\, i\le k\}.
\end{eqnarray*}
One can first show that the action of $\mathfrak{c}_{JG}(k,N)$ on $V$ is trivial. In fact, if this is not true,
since $\mathfrak{c}_{JG}(k,N)$ is an ideal of ${\frak g}(k,N)$,
then we have $V=\mathfrak{c}_{JG}(k,N)V$ by the irreducibility of $V$. Consider further the decomposition of ${\frak c}_{JG}(k,N)$:
$$
{\frak c}_{JG}(k,N)=\C \bar{J}_{k,N+1}+\check{{\frak c}}_{JG}(k,N), \mbox{\ \ where\ \ }\check{{\frak c}}_{JG}(k,N)={\frak c}_{JG}(k,N)\backslash\C \bar{J}_{k,N+1}.
$$
As in (2), by the comparing traces technique, we must have that the action of $\bar{J}_{k,N+1}$ is trivial, and so $V=\check{{\frak c}}_{JG}(k,N)V$.
Considering the action of $\bar{L}_{0,0}$ on $\check{{\frak c}}_{JG}(k,N)$, we still have \eqref{equ-L0J-c} and \eqref{equ-L0G-c}.
Since $\check{{\frak c}}_{JG}(k,N)$ is abelian, $\check{{\frak c}}_{JG}(k,N)$
acts nilpotently on $V$, which contradicts to $V=\check{{\frak c}}_{JG}(k,N)V$.

Now, by the above assumption, we have $V={\frak c}_{L}(k,N)V$.
Similarly, consider the decomposition of ${\frak c}_{L}(k,N)$:
$$
{\frak c}_{L}(k,N)=\C \bar{L}_{k,N}+\check{{\frak c}}_{L}(k,N), \mbox{\ \ where\ \ }\check{{\frak c}}_{L}(k,N)={\frak c}_{L}(k,N)\backslash\C \bar{L}_{k,N}.
$$
Note that $\bar{L}_{k,N}$ is an even central element of $\mathfrak{g}(k,N)$.
Here keep in mind that the actions of $\bar{J}_{k,N+1},\,\bar{G}_{k,N+\frac{1}{2}}^\pm\in\mathfrak{c}_{JG}(k,N)$ on $V$ are trivial.
Consider the actions of both sides of the following equation on $V$:
\begin{equation*}
[\bar{L}_{k,0}, \bar{L}_{0,N}] = -((N+1)k+Np)\bar{L}_{k,N}.
\end{equation*}
Again, by the comparing traces technique, one can show that the action of $\bar{L}_{k,N}$ is trivial.
Hence, $V=\check{{\frak c}}_L(k,N)V$.
Considering the action of $\bar{L}_{0,0}$ on $\check{{\frak c}}_L(k,N)$, we still have \eqref{equ-L0L-c}.
The only difference is the reason for $b_4\ne0$, which should be replaced by $i<k=j_L$.
Hence, all elements in $\check{{\frak c}}_L(k,N)$ act nilpotently on $V$. Since $\check{{\frak c}}_L(k,N)$ is abelian, $\check{{\frak c}}_L(k,N)$
acts nilpotently on $V$, which contradicts to $V=\check{{\frak c}}_L(k,N)V$.
\QED

\begin{lemm}\label{lemma-3}
Suppose $k,\,N\ge1$ and $\Omega(G)\ne \emptyset$.
Use the same notations $j_L,\,j_J,\,n_L,\,n_J$ as in Lemma~\ref{lemma-2}.
Let further
$$
j_G={\rm{max}}\{j\,|\,(j,n)\in\Omega(G)\}, \quad n_G={\rm{max}}\{n\,|\,(j,n)\in\Omega(G)\}.
$$
\baselineskip1pt\lineskip7pt\parskip-1pt
\begin{itemize}\parskip-1pt
  \item[{\rm(1)}] If $\max\{j_G, j_J\}<k$, then the row ideal $\mathfrak{r}(k,N)$ acts trivially on $V$;
  \item[{\rm(2)}] If $j_G<j_J=k$, $j_L<k$, then the row ideal $\mathfrak{r}(k,N)$ acts trivially on $V$;
  \item[{\rm(3)}] If $j_G<j_J=j_L=k$, $n_L<N$, then the column ideal $\mathfrak{c}(k,N)$ acts trivially on $V$;
  \item[{\rm(4)}] If $j_G<j_J=j_L=k$, $n_L=N$, then the column ideal $\mathfrak{c}(k,N)$ acts trivially on $V$;
  \item[{\rm(5)}] If $j_J<j_G=k$, $n_G<N$, then the column ideal $\mathfrak{c}(k,N)$ acts trivially on $V$;
  \item[{\rm(6)}] If $j_J<j_G=k$, $n_G=N$, then the row ideal $\mathfrak{r}(k,N)$ acts trivially on $V$.
\end{itemize}
\end{lemm}

\ni{\it Proof.}\ \
Note first that $j_J\ne j_G$, and so (1)--(6) cover all possible cases.

The conclusions (1)--(4) can be respectively proved in a similar way as Lemma~\ref{lemma-2}(1)--(4). The differences lie in the
reasons for $b_3\ne0$ and $b_6\ne0$, which should be replaced by $k>j_G$.

For (5), we have $j_L=j_J\, (<k)$, $n_L=n_G\, (<N)$, and $n_J=n_G-1$.
Assume that $\mathfrak{c}(k,N)$ acts non-trivially on $V$. We have
$V=\mathfrak{c}(k,N)V$ by the irreducibility of $V$.
Considering the action of $\bar{L}_{0,0}$ on $\mathfrak{c}(k,N)$, we still have
\eqref{equ-L0L-c}--\eqref{equ-L0G-c}. The only difference is the reason for $b_6\ne0$, which should be replaced by $N>n_G$.
The remaining arguments are the same as in Lemma~\ref{lemma-2}(3).

For (6), we have $k=p(N+\frac{1}{2})>pN\ge 2N$ (recall that $p\in 2\Z_{\ge1}$ by Proposition~\ref{prop-fact}(2)).
Assume that $\mathfrak{r}(k,N)$ acts non-trivially on $V$. Consider the decomposition of ${\frak r}(k,N)$:
\begin{eqnarray*}
{\frak r}(k,N) &\!\!\!=\!\!\!& \sp_{\C}\{\bar{J}_{k,N+1},\bar{G}_{k,N+\frac{1}{2}}^\pm\}+\tilde{{\frak r}}(k,N), \quad \mbox{where} \\
\tilde{{\frak r}}(k,N) &\!\!\!=\!\!\!& {\frak r}(k,N)\backslash\sp_{\C}\{\bar{J}_{k,N+1},\bar{G}_{k,N+\frac{1}{2}}^\pm\}.
\end{eqnarray*}
First, recall that $\bar{J}_{k,N+1}$ is an even central element of $\mathfrak{g}(k,N)$.
As in Lemma~\ref{lemma-2}(2), we must have that the action of $\bar{J}_{k,N+1}$ is trivial.
Then, by relations
\begin{eqnarray*}
[\bar{G}_{k,N+\frac{1}{2}}^{+}, \bar{G}_{0,\frac{1}{2}}^{-}] = -(N+1)p\bar{J}_{k,N+1}, \quad
[\bar{G}_{k,N+\frac{1}{2}}^{-}, \bar{G}_{0,\frac{1}{2}}^{+}] = (N+1)p\bar{J}_{k,N+1}.
\end{eqnarray*}
One can view $\bar{G}_{k,N+\frac{1}{2}}^\pm$ as odd central elements of $\mathfrak{g}(k,N)$.
We claim that the actions of $\bar{G}_{k,N+\frac{1}{2}}^\pm$ are also trivial.
If this is not true, then by Schur's lemma, we have $\dim V_{\bar{0}}=\dim V_{\bar{1}}$ and
the matrix of the action of $\bar{G}_{k,N+\frac{1}{2}}^+$ 
on $V$ with respect to a suitable basis has form
$$
\left(
  \begin{array}{cc}
     & -cI_d \\
    cI_d &  \\
  \end{array}
\right),
$$
where $c\ne0$, $d=\dim V_{\bar{0}}=\dim V_{\bar{1}}$, and $I_d$ denotes the identity matrix of order $d$.
Note further that the actions of $\bar{L}_{k,0}$ and $\bar{G}_{0,N+\frac{1}{2}}^+$ are nilpotent,
since the adjoint actions of $\bar{L}_{0,0}$ on them are non-trivial. We may assume that
the matrices of the actions of $\bar{L}_{k,0}$ and $\bar{G}_{0,N+\frac{1}{2}}^+$
on $V$ with respect to the above basis (adjust the even and odd parts of basis if necessary) respectively have forms
$$
\left(
  \begin{array}{cc}
    A_1 &  \\
      & A_2 \\
  \end{array}
\right) \quad\quad\mbox{and}\quad\quad
\left(
  \begin{array}{cc}
     & B_1 \\
    B_2 &  \\
  \end{array}
\right),
$$
where $A_1,\,A_2,\,B_1,\,B_2$ are matrices of order $d$, and $A_1$ (resp.~$A_2$) is a strictly upper (resp.~lower) triangular matrix.
Considering the actions of both sides of the following equation on $V$:
$$
[\bar{L}_{k,0}, \bar{G}_{0,N+\frac{1}{2}}^+]=b \bar{G}_{k,N+\frac{1}{2}}^+, \mbox{\ \ where\ \ } b=-(N+\frac{1}{2})(N+2)p\ne 0,
$$
we obtain
\begin{equation}\label{G-case6}
A_1B_1-B_1A_2=-bc I_d, \quad A_2B_2-B_2A_1=bc I_d.
\end{equation}
Comparing the entry $(d,d)$ of both sides of the first equation in \eqref{G-case6},
we obtain $0=-bc$, a contradiction.
Similarly, one can also prove the triviality of the action of $\bar{G}_{k,N+\frac{1}{2}}^-$. Thus, the above claim holds.
Then we have $V=\tilde{{\frak r}}(k,N)V$.
Furthermore, since $\tilde{{\frak r}}(k,N)$ is a completely reducible $\C \bar{L}_{0,0}$-module with no trivial summand
and $\tilde{{\frak r}}(k,N)$ is abelian, $\tilde{{\frak r}}(k,N)$ acts nilpotently on $V$, which contradicts to $V=\tilde{{\frak r}}(k,N)V$.
\QED
\vskip8pt

Now, we can give the proof of Theorem \ref{thm-1}.
\vskip8pt

\ni{\it Proof of Theorem \ref{thm-1}.}\ \
If $k=0$, then ${\frak g}(k,N)={\frak g}(0,N)={\cal A}(\KB)_+/{\cal I}(0,N)$.
It is equivalent to consider the problem for $\KB/\KB_{\langle 1\rangle}=\KB_{[0]}\cong K_2$.
By results in \cite{CL}, we must have that $\dim V = (1|1)$ or $(2|2)$.

Next, we assume that $k\ge 1$ and $N=0$. Note that the notations $\Omega(J)$, $\Omega(G)$ (cf.~\eqref{aux-J} and \eqref{aux-G})
and ${\frak r}(k,0)$ (cf.~\eqref{r-ideal}) still make sense.
We claim that the action of ${\frak r}(k,0)$ on $V$ is trivial.

{\it Case 1}: $p\notin \Z_+$.
In this case, we have $\Omega(J)=\Omega(G)=\emptyset$. One can prove the claim as in Lemma~\ref{lemma-2}(1).

{\it Case 2}: $p\in 1+2\Z_+$.
In this case, we have $\Omega(G)=\emptyset$.
If $j_J< k$, then $\Omega(J)=\emptyset$. One can prove the claim as in Lemma~\ref{lemma-2}(1).
If $j_J= k$, then $\Omega(J)\ne\emptyset$. One can prove the claim as in Lemma~\ref{lemma-2}(2).

{\it Case 3}: $p\in 2\Z_{\ge1}$.
Note that $j_J\ne j_G$.
If $\max\{j_G, j_J\}<k$, then $\Omega(J)=\Omega(G)=\emptyset$. One can prove the claim as in Lemma~\ref{lemma-2}(1).
If $j_J=k$, then $\Omega(J)\ne\emptyset$. One can prove the claim as in Lemma~\ref{lemma-2}(2).
If $j_G=k$, then $\Omega(G)\ne\emptyset$. One can prove the claim as in Lemma~\ref{lemma-3}(6).

The claim implies that $V$ is simply an irreducible module over $\mathfrak{g}(k-1,0)$.
By induction on $k$, we must have that $\dim V= (1|1)$ or $(2|2)$.

At last, we assume that $k,\,N\ge 1$. Note that
if the row ideal $\mathfrak{r}(k,N)$ (resp., column ideal $\mathfrak{c}(k,N)$) of
$\mathfrak{g}(k,N)$ acts trivially on $V$, then $V$ can be viewed as an irreducible module over
$\mathfrak{g}(k-1,N)$ (resp., $\mathfrak{g}(k,N-1)$).
By simultaneous induction on $k$ and $N$, using Lemmas~\ref{lemma-1}--\ref{lemma-3}, we must have that $\dim V=(1|1)$ or $(2|2)$.
\QED

\section{\large{Free conformal $\KB$-modules of small rank}} 

In this section, we classify all the free conformal modules of ranks $(1+1)$ and $(2+2)$ over $\KB$
by using the classification result of those over $K_2$.
For reducible ones, we completely determine their composition factors.

\subsection{Finite conformal modules over $K_2$}

Let us first list some finite conformal modules over $K_2=\langle L, J, G^\pm\rangle$.
For the convenience of later use, we choose $\{L_0=pL,\,J_0=J,\,G_0^\pm=\sqrt{p}G^\pm\}$ as a $\C[\partial]$-basis of $K_2$ (cf.~\eqref{equ-K2}).
\baselineskip1pt\lineskip7pt\parskip-1pt
\begin{itemize}\parskip-1pt
\item[{\rm(1)}] Trivial $K_2$-modules

For any $\a\in\C$, the even one-dimensional trivial $K_2$-module $\C c_\a$ satisfy $L_0\,{}_\lambda c_\a=J_0\,{}_\lambda c_\a=G^\pm_0\,{}_\lambda c_\a=0$
and $\partial c_\a=\a c_\a$, where $|c_\a|=\bar{0}$.

\item[{\rm(2)}] $K_2$-modules of rank $(1+1)$

For any $\D,\,\a\in\C$, there is a conformal $K_2$-module $K_{\D,\a}^{(1)}$ of rank $(1+1)$ with $\C[\partial]$-basis $\{v_{\bar{0}}, v_{\bar{1}}\}$
and $\lambda$-actions
\begin{equation}\label{K2-module-rk1+1-1}
\left\{\begin{array}{ll}
L_0\,{}_\lambda\, v_{\bar{0}}=p(\partial+\D \lambda+\a)v_{\bar{0}}, &  L_0\,{}_\lambda\, v_{\bar{1}}=p(\partial+(\D+\frac{1}{2}) \lambda+\a)v_{\bar{1}},  \\[3pt]
J_0\,{}_\lambda\, v_{\bar{0}}=-2 \D v_{\bar{0}}, &  J_0\,{}_\lambda\, v_{\bar{1}}=(1-2\D) v_{\bar{1}},  \\[3pt]
G_0^+\,{}_\lambda\, v_{\bar{0}}=\sqrt{p}v_{\bar{1}}, & G_0^+\,{}_\lambda\, v_{\bar{1}}=0, \\[3pt]
G_0^-\,{}_\lambda\, v_{\bar{0}}=0, & G_0^-\,{}_\lambda\, v_{\bar{1}}=2\sqrt{p}(\partial+2\D\lambda+\a)v_{\bar{0}}.
\end{array}\right.
\end{equation}
There is another conformal $K_2$-module $K_{\D,\a}^{(2)}$ of rank $(1+1)$ with $\C[\partial]$-basis $\{v_{\bar{0}}, v_{\bar{1}}\}$
and $\lambda$-actions
\begin{equation}\label{K2-module-rk1+1-2}
\left\{\begin{array}{ll}
L_0\,{}_\lambda\, v_{\bar{0}}=p(\partial+\D \lambda+\a)v_{\bar{0}}, &  L_0\,{}_\lambda\, v_{\bar{1}}=p(\partial+(\D+\frac{1}{2}) \lambda+\a)v_{\bar{1}},  \\[3pt]
J_0\,{}_\lambda\, v_{\bar{0}}=2 \D v_{\bar{0}}, &  J_0\,{}_\lambda\, v_{\bar{1}}=(2\D-1) v_{\bar{1}},  \\[3pt]
G_0^+\,{}_\lambda\, v_{\bar{0}}=0, & G_0^+\,{}_\lambda\, v_{\bar{1}}=2\sqrt{p}(\partial+2\D\lambda+\a)v_{\bar{0}}, \\[3pt]
G_0^-\,{}_\lambda\, v_{\bar{0}}=\sqrt{p}v_{\bar{1}}, & G_0^-\,{}_\lambda\, v_{\bar{1}}=0.
\end{array}\right.
\end{equation}

\item[{\rm(3)}] $K_2$-modules of rank $(2+2)$

For any $\D,\,\L,\,\a\in\C$, there is a conformal $K_2$-module
$K_{\D,\L,\a}$ of rank $(2+2)$ with $\C[\partial]$-basis $\{v_{\bar{0}}^{(1)}, v_{\bar{0}}^{(2)}, v_{\bar{1}}^{(1)}, v_{\bar{1}}^{(2)}\}$
and $\lambda$-actions
\begin{equation}\label{K2-module-rk2+2}
\left\{\begin{array}{ll}
L_0\,{}_\lambda\, v_{\bar{0}}^{(1)}=p(\partial+\D \lambda+\a)v_{\bar{0}}^{(1)},
\qquad\   L_0\,{}_\lambda\, v_{\bar{1}}^{(\ell)}=p(\partial+(\D+\frac{1}{2}) \lambda+\a)v_{\bar{1}}^{(\ell)},\ \ell=1,2,  \\[4pt]
L_0\,{}_\lambda\, v_{\bar{0}}^{(2)}=p(\partial+(\D+1) \lambda+\a)v_{\bar{0}}^{(2)}+p(\D+\frac{\L}{2})\lambda^2 v_{\bar{0}}^{(1)},  \\[4pt]
J_0\,{}_\lambda\, v_{\bar{0}}^{(1)}=\L v_{\bar{0}}^{(1)}, \qquad\qquad\qquad\qquad\   J_0\,{}_\lambda\, v_{\bar{1}}^{(1)}=(\L+1) v_{\bar{1}}^{(1)},  \\[4pt]
J_0\,{}_\lambda\, v_{\bar{0}}^{(2)}=\L v_{\bar{0}}^{(2)}+(2\D+\L)\lambda v_{\bar{0}}^{(1)}, \quad   J_0\,{}_\lambda\, v_{\bar{1}}^{(2)}=(\L-1) v_{\bar{1}}^{(2)},  \\[4pt]
G_0^+\,{}_\lambda\, v_{\bar{0}}^{(1)}=\sqrt{p}v_{\bar{1}}^{(1)}, \qquad\qquad\qquad\quad G_0^+\,{}_\lambda\, v_{\bar{1}}^{(1)}=0, \\[4pt]
G_0^+\,{}_\lambda\, v_{\bar{0}}^{(2)}=-\sqrt{p}(2\D+\L)\lambda v_{\bar{1}}^{(1)},
\quad\ \  G_0^+\,{}_\lambda\, v_{\bar{1}}^{(2)}=\sqrt{p}(2\D+\L)\lambda v_{\bar{0}}^{(1)} + \sqrt{p} v_{\bar{0}}^{(2)}, \\[4pt]
G_0^-\,{}_\lambda\, v_{\bar{0}}^{(1)}=\sqrt{p}v_{\bar{1}}^{(2)},
\qquad\qquad\qquad\quad G_0^-{}_\lambda\, v_{\bar{1}}^{(1)}\!=\!\sqrt{p}(2\partial\!+\!(2\D\!-\!\L)\lambda\!+\!2\a)v_{\bar{0}}^{(1)}\!\!-\!\!\sqrt{p}v_{\bar{0}}^{(2)},\\[4pt]
G_0^-\,{}_\lambda\, v_{\bar{0}}^{(2)}=\sqrt{p}(2\partial\!+\!(2\D\!-\!\L\!+\!2)\lambda\!+\!2\a)v_{\bar{1}}^{(2)},
\ \, G_0^-\,{}_\lambda\, v_{\bar{1}}^{(2)}=0.
\end{array}\right.
\end{equation}
\end{itemize}

\begin{lemm}\label{lemma-FICMs-K2}
Let $V$ be a finite irreducible conformal module over $K_2$. Up to parity change, $V$ is isomorphic to
$\C c_\a$, $K_{\D,\a}^{(1)}$, $K_{\D,\a}^{(2)}$ with $\D\ne0$, or $K_{\D,\L,\a}$ with $2\D\pm\L\ne0$ listed above.
\end{lemm}

\begin{rema}\label{remark-K2}
\rm
The above classification result of FICMs over $K_2$ was given in \cite{CL}.
By using a lemma in \cite{CK} (see Lemma~\ref{lemma-general}) or the  arguments in \cite[Remark~8.3]{DK}, we have that
\baselineskip1pt\lineskip7pt\parskip-1pt
\begin{itemize}\parskip-1pt
\item[{\rm(1)}] up to parity change, any conformal $K_2$-module of rank $(1+1)$ has form $K_{\D,\a}^{(1)}$ or $K_{\D,\a}^{(2)}$;
\item[{\rm(2)}] up to parity change, any conformal $K_2$-module of rank $(2+2)$ has form $K_{\D,\L,\a}$.
\end{itemize}
Furthermore, $K_{\D,\a}^{(1)}$ and $K_{\D,\a}^{(2)}$ are irreducible if and only if $\D\ne0$,
while $K_{\D,\L,\a}$ is irreducible if and only if $2\D\pm\L\ne0$.
The composition factors of reducible $K_2$-modules $K_{0,\a}^{(1)}$, $K_{0,\a}^{(2)}$, $K_{\D,\pm 2\D,\a}$ will be given in Remark~\ref{remark-K2-CF}.
\end{rema}

\subsection{$\KB$-modules of rank $(1+1)$}

In this subsection, we classify all the free conformal $\KB$-modules of rank $(1+1)$.
Obviously, for any $\D,\,\a\in\C$, the conformal $K_2$-modules $K_{\D,\a}^{(1)}$ and $K_{\D,\a}^{(2)}$
can be trivially extended to conformal $\KB$-modules $V_{\D,\a}^{(1)}$ and $V_{\D,\a}^{(2)}$:
\baselineskip1pt\lineskip7pt\parskip-1pt
\begin{itemize}\parskip-1pt
  \item[{\rm(I-1)}] $V_{\D,\a}^{(1)}=\C[\partial]v_{\bar{0}}\oplus\C[\partial]v_{\bar{1}}$ with \eqref{K2-module-rk1+1-1} and
  $L_i\,{}_\lambda\, v_{s}= J_i\,{}_\lambda\, v_{s}= G_i^\pm\,{}_\lambda\,v_{s}=0,\, i\ge 1,\, s\in\Z/2\Z$;
  \item[{\rm(I-2)}] $V_{\D,\a}^{(2)}=\C[\partial]v_{\bar{0}}\oplus\C[\partial]v_{\bar{1}}$ with \eqref{K2-module-rk1+1-2} and
  $L_i\,{}_\lambda\, v_{s}= J_i\,{}_\lambda\, v_{s}= G_i^\pm\,{}_\lambda\,v_{s}=0,\, i\ge 1,\, s\in\Z/2\Z$.
\end{itemize}

\begin{lemm}\label{thm-2}
Let $M$ be a nontrivial finite conformal module over $\KB$. Then the $\lambda$-actions of $L_i, J_i, G_i^\pm\in\KB$ on $M$ are trivial for all $i\gg 0$.
In particular, a finite conformal module over $\KB$ is simply a finite conformal module over $\KB_{[n]}$ for some big enough integer $n$,
where $\KB_{[n]}$ is defined by \eqref{equ-quotiont}.
\end{lemm}

\ni{\it Proof.}\ \
Since $M$ can be viewed as a conformal module over the conformal subalgebra $\langle L_i\,|\,i\in\Z_+\rangle$ of $\KB$, it follows from \cite{SXY} that
the $\lambda$-action of $L_i$ on $M$ is trivial for all $i\gg 0$.
Furthermore, from \eqref{brackets-LJ} and \eqref{brackets-LG} we see that the $\lambda$-actions of $J_i, G_i^\pm$ on $M$ are also trivial for all $i\gg 0$.
(Here, we would like to take this opportunity to point out that \cite[Theorem~4.2]{X} for $\SB$ can be also concisely proved by using the above observation.)
\QED

\begin{theo}\label{thm-rank-1+1}
Let $M$ be a nontrivial free conformal module of rank $(1+1)$ over $\KB$.
Then, up to parity change, $M\cong V_{\D,\a}^{(1)}$ or $V_{\D,\a}^{(2)}$ for some $\D,\a\in\C$.
\end{theo}

\ni{\it Proof.}\ \
Let $M=\C[\partial]v_{\bar{0}}\oplus\C[\partial]v_{\bar{1}}$. By regarding $M$ as a conformal module over $K_2$,
we see that, up to parity change, the $\lambda$-actions of $L_0, J_0, G_0^\pm$
have forms \eqref{K2-module-rk1+1-1} or \eqref{K2-module-rk1+1-2}, where $\D,\a\in\C$.
By Lemma~\ref{thm-2}, $L_i\,{}_\lambda\, v_s=J_i\,{}_\lambda\, v_s=G_i^\pm\,{}_\lambda\, v_s=0$ for $i\gg 0$, $s\in\Z/2\Z$.
Note that $\KB$ is $\Z$-graded in the sense that $\KB=\oplus_{i\in \Z_+}\KB_i$, where $\KB_i=\C[\partial]L_i\oplus\C[\partial]J_i\oplus\C[\partial]G_i^+\oplus\C[\partial]G_i^-$.
Assume that $k\in\Z_+$ is the largest integer such that the action of $\KB_k$ on $M$ is nontrivial.

If $k=0$, then $M$ is simply a conformal $K_2$-module. By Remark~\ref{remark-K2}(1), the conclusion holds.

Next, we consider the case $k>0$. By assumption, we can suppose
$$
L_k\,{}_\lambda\, v_{s}=a_{s}(\partial,\lambda)v_{s}, \quad J_k\,{}_\lambda\, v_{s}=b_{s}(\partial,\lambda)v_{s}, \quad
G_k^\pm\,{}_\lambda\, v_{s}=c_{s}^\pm(\partial,\lambda)v_{s+\bar{1}},
$$
where $a_{s}(\partial,\lambda),\,b_{s}(\partial,\lambda),\,c_{s}^\pm(\partial,\lambda)\in\C[\partial, \lambda]$ for $s\in\Z/2\Z$,
and at least one of them is nonzero.
Considering the action of $[L_k\,{}_\lambda\, L_k]\,{}_{\lambda+\mu}=0$ on $v_{s}$, we obtain
$$
a_{s}(\partial,\lambda)a_{s}(\partial+\lambda,\mu) = a_{s}(\partial,\mu)a_{s}(\partial+\mu,\lambda),\ \ s\in\Z/2\Z.
$$
Comparing the coefficients of $\lambda$,
we see that $a_{s}(\partial,\lambda)$ is independent of the variable $\partial$,
and so we can denote $a_{s}(\lambda)=a_{s}(\partial,\lambda)$.
Then, considering the action of the operator $[L_0\,{}_\lambda\, L_k]\,{}_{\lambda+\mu}=((k+p)\lambda-p\mu)L_{k}\,{}_{\lambda+\mu}$
on $v_{s}$, we obtain
\begin{equation}
\label{equ-a-type} (p\mu-(k+p)\lambda)a_{s}(\lambda+\mu) = p\mu a_{s}(\mu),\ \ s\in\Z/2\Z.
\end{equation}

If $k\ne -p$, then $k+p\ne 0$. By \eqref{equ-a-type} with $\mu=0$, we obtain $a_{s}(\lambda)=0$.
Hence, the action of $L_k$ on $M$ is trivial. Furthermore, by relations
\begin{eqnarray}
\nonumber (k+p)(\partial+\lambda)J_k &\!\!\!=\!\!\!& [L_k\,{}_\lambda\, J_0],\\
\nonumber ((k+p)\partial+(k+\frac{3}{2}p)\lambda)G_k^\pm &\!\!\!=\!\!\!& [L_k\,{}_\lambda\, G_0^\pm],
\end{eqnarray}
we see that the actions of $J_k$ and $G_k^\pm$ on $M$ are also trivial, a contradiction.

If $k=-p$, then $p$ is a negative integer. By \eqref{equ-a-type},  we see that  $a_{s}(\lambda)$
is independent of the variable $\lambda$, and so we can denote $a_{s}=a_{s}(\lambda), s\in\Z/2\Z$.

If $p\le -2$, then $k\ge 2$.
We claim that the action of $\KB_{k-1}$ on $M$ is trivial. As in \cite[Thenrem~5.1]{X}, one can first prove that the action of $L_{k-1}$ on $M$ is trivial.
Then, by relations
\begin{eqnarray}
\nonumber -(\partial+\lambda)J_{k-1} &\!\!\!=\!\!\!& [L_{k-1}\,{}_\lambda\, J_0],\\
\nonumber -(\partial+(1-\frac{1}{2}p)\lambda)G_{k-1}^\pm &\!\!\!=\!\!\!& [L_{k-1}\,{}_\lambda\, G_0^\pm],
\end{eqnarray}
we see that the actions of $J_{k-1}$ and $G_{k-1}^\pm$ on $M$ are also trivial, and thus the claim holds.
At last, note that in case $p\le -2$, the $\C[\partial]$-basis of $\KB_k$ can be generated by those of $\KB_{k-1}$ and $L_1$.
Hence, the action of $\KB_k$ on $M$ is trivial, a contradiction.

If $p=-1$, then $k=1$. If $M\cong K_{\D,\a}^{(1)}$ as a conformal $K_2$-module, then the actions of $L_0$, $J_0$ and $G_0^\pm$ have
the form of \eqref{K2-module-rk1+1-1}.
Applying the operator $[L_1\,{}_\lambda G_{0}^+]\,{}_{\lambda+\mu}=-\frac{1}{2}\lambda G_1^+\,{}_{\lambda+\mu}$ on $v_{\bar{0}}$, we obtain
$$
\sqrt{-1} (a_{\bar{0}}-a_{\bar{1}})v_{\bar{1}}=\frac{1}{2}\lambda\, c_{\bar{0}}^+(\partial, \lambda+\mu)v_{\bar{1}},
$$
which implies that $a_{\bar{0}}=a_{\bar{1}}$ and $c_{\bar{0}}^+(\partial, \lambda)=0$.
Then, applying the above operator on $v_{\bar{1}}$, we obtain
$$
0=-\frac{1}{2}\lambda\, c_{\bar{1}}^+(\partial, \lambda+\mu)v_{\bar{0}},
$$
which implies that $c_{\bar{1}}^+(\partial, \lambda)=0$. Denote $\b=a_{\bar{0}}=a_{\bar{1}}$.
Applying the operator $[L_1\,{}_\lambda G_{0}^-]\,{}_{\lambda+\mu}=-\frac{1}{2}\lambda G_1^-\,{}_{\lambda+\mu}$ on $v_{\bar{0}}$ and $v_{\bar{1}}$, respectively,
we obtain
$$
0=-\frac{1}{2}\lambda\, c_{\bar{0}}^-(\partial, \lambda+\mu)v_{\bar{1}},\quad 2\sqrt{-1}\b\lambda v_{\bar{0}}=-\frac{1}{2}\lambda\, c_{\bar{1}}^-(\partial, \lambda+\mu)v_{\bar{0}},
$$
which imply that $c_{\bar{0}}^-(\partial, \lambda)=0$ and $c_{\bar{1}}^-(\partial, \lambda)=-4\sqrt{-1}\b$.
Applying the operator $[G_0^+\,{}_\lambda G_{1}^-]\,{}_{\lambda+\mu}=(-\partial J_1+2L_1)\,{}_{\lambda+\mu}$ on $v_{s}$,
we obtain
$$
4\b=(\lambda+\mu)b_{s}(\partial, \lambda+\mu) +2\b, \ \ s\in\Z/2\Z,
$$
which imply that $\b=0$ and $b_{s}(\partial, \lambda)=0$ for $s\in\Z/2\Z$. Hence, $M\cong V_{\D,\a}^{(1)}$.
Similarly, if $M\cong K_{\D,\a}^{(2)}$ as a conformal $K_2$-module, then one can show that $M\cong V_{\D,\a}^{(2)}$.
This completes the proof.
\QED

\begin{rema}\label{remark-1+1-triviality}
\rm
In \cite{X}, we showed that a rank $(1+1)$ $K_1$-module can be nontrivially extended to a $\SB$-module if $p=-1$.
However, by Theorem~\ref{thm-rank-1+1}, we see that any rank $(1+1)$ $K_2$-module $M$
can not be nontrivially extended to a $\KB$-module, even in case that $p=-1$ and $M$ can be degenerated from a rank $(2+2)$ $K_2$-module
(note that a $K_2$-module $K_{\D,\L,\a}$ will degenerate to a $K_2$-module $K_{\D,\a}^{(1)}$
if we set $v_{\bar{0}}^{(2)}=v_{\bar{1}}^{(2)}=0$ and $\L=-2\D$). We shall see in the next subsection that
a rank $(2+2)$ $K_2$-module can be nontrivially extended to a $\KB$-module if $p=-1$.
\end{rema}

\subsection{$\KB$-modules of rank $(2+2)$}

In this subsection, we classify all the free conformal $\KB$-modules of rank $(2+2)$.
As above, for any $\D,\,\L,\,\a\in\C$, the conformal $K_2$-module $K_{\D,\L,\a}$ can be trivially extended to a conformal $\KB$-module $V_{\D,\L,\a}$,
while for $\KBB(-1)$, the conformal $K_2$-module $K_{\D,\L,\a}$ can be extended to a conformal $\KB$-module $V_{\D,\L,\a,\b}$,
which is a nontrivial extension of $K_{\D,\L,\a}$ if $\b\ne 0$:
\baselineskip1pt\lineskip7pt\parskip-1pt
\begin{itemize}\parskip-1pt
  \item[{\rm(II-1)}] $V_{\D,\L,\a}=\C[\partial]v_{\bar{0}}^{(1)}\oplus\C[\partial]v_{\bar{0}}^{(2)}\oplus\C[\partial]v_{\bar{1}}^{(1)}\oplus\C[\partial]v_{\bar{1}}^{(2)}$
  with \eqref{K2-module-rk2+2} and
  $L_i\,{}_\lambda\, v_{s}^{(\ell)}= J_i\,{}_\lambda\, v_{s}^{(\ell)}= G_i^\pm\,{}_\lambda\,v_{s}^{(\ell)}=0,\, i\ge 1,\, s\in\Z/2\Z,\, \ell=1,2$;
  \item[{\rm(II-2)}] $V_{\D,\L,\a,\b}=\C[\partial]v_{\bar{0}}^{(1)}\oplus\C[\partial]v_{\bar{0}}^{(2)}\oplus\C[\partial]v_{\bar{1}}^{(1)}\oplus\C[\partial]v_{\bar{1}}^{(2)}$
   with \eqref{K2-module-rk2+2} and
\begin{equation}\label{KB-non-trivial-extension}
\left\{\begin{array}{ll}
   L_1\,{}_\lambda\, v_{\bar{0}}^{(1)}=-\frac{\sqrt{-1}}{2}\b v_{\bar{0}}^{(1)},
&  L_1\,{}_\lambda\, v_{\bar{1}}^{(1)}=-\frac{\sqrt{-1}}{2}\b v_{\bar{1}}^{(1)},  \\[4pt]
   L_1\,{}_\lambda\, v_{\bar{0}}^{(2)}=-\frac{\sqrt{-1}}{2}\b (v_{\bar{0}}^{(2)}+\lambda v_{\bar{0}}^{(1)}),
&  L_1\,{}_\lambda\, v_{\bar{1}}^{(2)}=-\frac{\sqrt{-1}}{2}\b v_{\bar{1}}^{(2)},  \\[4pt]
   J_1\,{}_\lambda\, v_{\bar{0}}^{(1)}=0, &  J_1\,{}_\lambda\, v_{\bar{1}}^{(1)}=0,  \\[4pt]
   J_1\,{}_\lambda\, v_{\bar{0}}^{(2)}=\sqrt{-1}\b v_{\bar{0}}^{(1)}, &  J_1\,{}_\lambda\, v_{\bar{1}}^{(2)}=0,  \\[4pt]
   G_1^+\,{}_\lambda\, v_{\bar{0}}^{(1)}=0, & G_1^+\,{}_\lambda\, v_{\bar{1}}^{(1)}=0, \\[4pt]
   G_1^+\,{}_\lambda\, v_{\bar{0}}^{(2)}=\b v_{\bar{1}}^{(1)}, & G_1^+\,{}_\lambda\, v_{\bar{1}}^{(2)}=-\b v_{\bar{0}}^{(1)}, \\[4pt]
   G_1^-\,{}_\lambda\, v_{\bar{0}}^{(1)}=0, & G_1^-\,{}_\lambda\, v_{\bar{1}}^{(1)}=-\b v_{\bar{0}}^{(1)}, \\[4pt]
   G_1^-\,{}_\lambda\, v_{\bar{0}}^{(2)}=-\b v_{\bar{1}}^{(2)}, & G_1^-\,{}_\lambda\, v_{\bar{1}}^{(2)}=0, \\[4pt]
   X_i\,{}_\lambda\, v_{s}^{(\ell)}=0,
\end{array}\right.
\end{equation}
where $X=L,J,G^\pm,\, i\ge 2,\, s\in\Z/2\Z,\, \ell=1,2$.
\end{itemize}

\begin{theo}\label{thm-rank-2+2}
Let $M$ be a nontrivial free conformal module of rank $(2+2)$ over $\KB$.
\baselineskip1pt\lineskip7pt\parskip-1pt
\begin{itemize}\parskip-1pt
  \item[{\rm(1)}] If $p\ne -1$, then, up to parity change, $M\cong V_{\D,\L,\a}$ for some $\D,\L,\a\in\C$.
  \item[{\rm(2)}] If $p=-1$, then, up to parity change, $M\cong V_{\D,\L,\a,\b}$ for some $\D,\L,\a,\b\in\C$.
\end{itemize}
\end{theo}

\ni{\it Proof.}\ \
Suppose $M$ is a nontrivial free conformal $\KB$-module of rank $(2+2)$.
Write $M=M_{\bar{0}}\oplus M_{\bar{1}}$, where $M_{s}=\C[\partial]v_{s}^{(1)}\oplus\C[\partial]v_{s}^{(2)}$, $s\in\Z/2\Z$.
By regarding $M$ as a conformal module over $K_2$,
we see that, up to parity change, the $\lambda$-actions of $L_0, J_0, G_0^\pm$
have forms \eqref{K2-module-rk2+2}, where $\D,\L,\a\in\C$.
By Lemma~\ref{thm-2}, $L_i\,{}_\lambda\, v_s=J_i\,{}_\lambda\, v_s=G_i^\pm\,{}_\lambda\, v_s=0$ for $i\gg 0$, $s\in\Z/2\Z$.
Assume that $k\in\Z_+$ is the largest integer such that the action of $\KB_k$ on $M$ is nontrivial.

If $k=0$, then $M$ is simply a conformal $K_2$-module. By Remark~\ref{remark-K2}(2), Theorem~\ref{thm-rank-2+2} holds.
More precisely, up to parity change,
\begin{equation*}
M\cong
\left\{\begin{array}{ll}
V_{\D,\L,\a}, & \mbox{if} \quad p\ne -1;\\[4pt]
V_{\D,\L,\a,0}, & \mbox{if} \quad p= -1.
\end{array}\right.
\end{equation*}

Next, we consider the case $k>0$.
By assumption, we can suppose
$$
\begin{aligned}
L_k\,{}_\lambda\, v_{s}^{(\ell)} &= a_{s}^{(\ell)}(\partial,\lambda)v_{s}^{(1)} + d_{s}^{(\ell)}(\partial,\lambda)v_{s}^{(2)}, \\[-2pt]
J_k\,{}_\lambda\, v_{s}^{(\ell)} &= b_{s}^{(\ell)}(\partial,\lambda)v_{s}^{(1)} + e_{s}^{(\ell)}(\partial,\lambda)v_{s}^{(2)}, \\[-2pt]
G_k^\pm\,{}_\lambda\, v_{s}^{(\ell)} &= c_{s}^{\pm(\ell)}(\partial,\lambda)v_{s+\bar{1}}^{(1)} + f_{s}^{\pm(\ell)}(\partial,\lambda)v_{s+\bar{1}}^{(2)},
\end{aligned}
$$
where $x_{s}^{(\ell)}(\partial,\lambda)\in\C[\partial, \lambda]$ for $x=a,b,c^\pm,d,e,f^\pm$, $s\in\Z/2\Z$, $\ell=1,2$,
and at least one of them is nonzero.
Let us first determine the actions of $G_k^+$ and $J_k$ on $M$ through three lemmas.

\begin{lemm}\label{action-G-on-M}
There exist $\b,\gamma\in\C$ such that
\baselineskip1pt\lineskip7pt\parskip-1pt
\begin{itemize}\parskip-1pt
  \item[{\rm(1)}] $G_k^+\,{}_\lambda\, v_{\bar{1}}^{(1)}=0$, $G_k^+\,{}_\lambda\, v_{\bar{0}}^{(1)}=\delta^*_\lambda v_{\bar{1}}^{(1)}$,
  where $\delta^*_\lambda=\delta_{k,p}\delta_{\D,\frac{3}{2}\L+1}\gamma(\partial+(\L+1)\lambda+\a)$;
  \item[{\rm(2)}] $G_k^+\,{}_\lambda\, v_{\bar{0}}^{(2)}=\delta_{k+p,0}\b v_{\bar{1}}^{(1)}$,
  $G_k^+\,{}_\lambda\, v_{\bar{1}}^{(2)}=-\delta_{k+p,0}\b v_{\bar{0}}^{(1)}$.
\end{itemize}
\end{lemm}

\ni{\it Proof.}\ \
(1) Considering the action of $[J_0\,{}_\lambda\, G_k^+]\,{}_{\lambda+\mu}=G_k^+{}_{\lambda+\mu}$ on $v_{\bar{0}}^{(1)}$, we obtain
\begin{eqnarray}
\label{J0Gk-V01c} (\L+1)c_{\bar{0}}^{+(1)}(\partial+\lambda,\mu)- \L c_{\bar{0}}^{+(1)}(\partial,\mu) &\!\!\!=\!\!\!& c_{\bar{0}}^{+(1)}(\partial,\lambda+\mu), \\
\label{J0Gk-V01f} (\L-1)f_{\bar{0}}^{+(1)}(\partial+\lambda,\mu)- \L f_{\bar{0}}^{+(1)}(\partial,\mu) &\!\!\!=\!\!\!& f_{\bar{0}}^{+(1)}(\partial,\lambda+\mu).
\end{eqnarray}
Taking $\lambda=0$ in \eqref{J0Gk-V01f}, we obtain $f_{\bar{0}}^{+(1)}(\partial,\mu)=0$.
Then, considering the action of $[G_0^+\,{}_\lambda\, G_k^+]\,{}_{\lambda+\mu}=0$ on $v_{\bar{0}}^{(1)}$, we obtain
$c_{\bar{1}}^{+(1)}(\partial,\mu)=f_{\bar{1}}^{+(1)}(\partial,\mu)=0$. Namely, $G_k^+\,{}_\lambda\, v_{\bar{1}}^{(1)}=0$.
Furthermore, applying $[L_0\,{}_\lambda\, G_k^+]\,{}_{\lambda+\mu}=((k+\frac{1}{2}p)\lambda-p\mu)G_k^+{}_{\lambda+\mu}$ on $v_{\bar{0}}^{(1)}$, we obtain
$$
p(\partial+(\D+\frac{1}{2})\lambda+\a)c_{\bar{0}}^{+(1)}(\partial+\lambda,\mu)-p(\partial+\mu+\D\lambda+\a)c_{\bar{0}}^{+(1)}(\partial,\mu)
=((k+\frac{1}{2}p)\lambda-p\mu)c_{\bar{0}}^{+(1)}(\partial,\lambda+\mu).
$$
Taking $\mu=0$, we obtain
\begin{equation}\label{L0Gk-V01}
p(\partial+\D\lambda+\a)\frac{c_{\bar{0}}^{+(1)}(\partial+\lambda,0)-c_{\bar{0}}^{+(1)}(\partial,0)}{\lambda}
+\frac{p}{2} c_{\bar{0}}^{+(1)}(\partial+\lambda,0)
=(k+\frac{1}{2}p)c_{\bar{0}}^{+(1)}(\partial,\lambda).
\end{equation}
Taking $\lambda\rightarrow0$, we obtain
\begin{equation}\label{L0Gk-V01-limit}
p(\partial+\a)\frac{d}{d\partial}c_{\bar{0}}^{+(1)}(\partial,0)=k c_{\bar{0}}^{+(1)}(\partial,0).
\end{equation}
Let $\deg\, c_{\bar{0}}^{+(1)}(\partial,0)=n$. The solution to \eqref{L0Gk-V01-limit} is $k=pn$ and $c_{\bar{0}}^{+(1)}(\partial,0)=\gamma_n(\partial+\a)^n$,
where $\gamma_n\in\C$ and $\gamma_n\ne0$ if $n\ge 1$.
By \eqref{J0Gk-V01c} with $\mu=0$ and \eqref{L0Gk-V01}, we have
$$
(\partial+\lambda+\a)^n\left(((n+\frac{1}{2})\L+n-\D)\lambda-\partial-\a\right)\gamma_n
=(\partial+\a)^n\left(((n+\frac{1}{2})\L-\D)\lambda-\partial-\a\right)\gamma_n.
$$
If $n\ge2$, by comparing the coefficients of $\lambda^{n+1}$ in the above equation, we must have $(n+\frac{1}{2})\L+n-\D=0$.
Then, by comparing the coefficients of $\lambda^{n}$, we obtain $-(\partial+\a)\gamma_n=0$, a contradiction.
If $n=1$, then $k=p$. By comparing the coefficients of $\lambda^2$ in the above equation, we have $\frac{3}{2}\L+1-\D=0$.
Then, by \eqref{J0Gk-V01c} with $\mu=0$, we obtain $c_{\bar{0}}^{+(1)}(\partial,\lambda)=\gamma_1(\partial+(\L+1)\lambda+\a)$.
If $n=0$, using \eqref{L0Gk-V01-limit} and then \eqref{J0Gk-V01c} with $\mu=0$, we obtain $c_{\bar{0}}^{+(1)}(\partial,\lambda)=0$.
Namely, we have
\begin{equation*}
G_k^+\,{}_\lambda\, v_{\bar{0}}^{(1)}=
\left\{\begin{array}{ll}
\gamma_1(\partial+(\L+1)\lambda+\a)v_{\bar{1}}^{(1)}, & \mbox{\rm if} \quad k=p,\, \D=\frac{3}{2}\L+1;\\[4pt]
0, & \mbox{\rm otherwise}.
\end{array}\right.
\end{equation*}
Denote $\delta^*_\lambda=\delta_{k,p}\delta_{\D,\frac{3}{2}\L+1}\gamma(\partial+(\L+1)\lambda+\a)$, where $\gamma=\gamma_1$.
Then, we have $G_k^+\,{}_\lambda\, v_{\bar{0}}^{(1)}=\delta^*_\lambda v_{\bar{1}}^{(1)}$.

(2) Applying  $[J_0\,{}_\lambda\, G_k^+]\,{}_{\lambda+\mu}=G_k^+{}_{\lambda+\mu}$ on $v_{\bar{0}}^{(2)}$, we obtain
\begin{eqnarray}
\label{J0Gk-V02c} (\L+1)c_{\bar{0}}^{+(2)}(\partial+\lambda,\mu)- \L c_{\bar{0}}^{+(2)}(\partial,\mu) &\!\!\!=\!\!\!& c_{\bar{0}}^{+(2)}(\partial,\lambda+\mu)
-(2\D+\L)\lambda\delta^*_{\mu}, \\
\label{J0Gk-V02f} (\L-1)f_{\bar{0}}^{+(2)}(\partial+\lambda,\mu)- \L f_{\bar{0}}^{+(2)}(\partial,\mu) &\!\!\!=\!\!\!& f_{\bar{0}}^{+(2)}(\partial,\lambda+\mu).
\end{eqnarray}
Taking $\lambda=0$ in \eqref{J0Gk-V02f}, we obtain $f_{\bar{0}}^{+(2)}(\partial,\mu)=0$.
Then, similar to \eqref{L0Gk-V01}, by considering the operator
$[L_0\,{}_\lambda\, G_k^+]\,{}_{\lambda+\mu}=((k+\frac{1}{2}p)\lambda-p\mu)G_k^+{}_{\lambda+\mu}$ on $v_{\bar{0}}^{(2)}$,
we obtain
\begin{eqnarray}
\nonumber
&\!\!\!\!\!\!& p(\partial+(\D+1)\lambda+\a)\frac{c_{\bar{0}}^{+(2)}(\partial+\lambda,0)-c_{\bar{0}}^{+(2)}(\partial,0)}{\lambda}
-\frac{p}{2} c_{\bar{0}}^{+(2)}(\partial+\lambda,0)\\
\label{L0Gk-V02} &\!\!\!=\!\!\!& (k+\frac{1}{2}p)c_{\bar{0}}^{+(2)}(\partial,\lambda)+p(\D+\frac{\L}{2})\lambda \delta^*_0.
\end{eqnarray}
Taking $\lambda\rightarrow0$, we obtain
\begin{equation}\label{L0Gk-V02-limit}
p(\partial+\a)\frac{d}{d\partial}c_{\bar{0}}^{+(2)}(\partial,0)=(k+p) c_{\bar{0}}^{+(2)}(\partial,0).
\end{equation}
Following the discussion in (1), by \eqref{J0Gk-V02c}, \eqref{L0Gk-V02} and \eqref{L0Gk-V02-limit}, one can show that
$c_{\bar{0}}^{+(2)}(\partial,\lambda)=\delta_{k+p,0}\b$ for some $\b\in\C$.
Namely, the action of $G_k^+$ on $v_{\bar{0}}^{(2)}$ has the required form.
Furthermore, considering the action of $[G_0^+\,{}_\lambda\, G_k^+]\,{}_{\lambda+\mu}=0$ on $v_{\bar{1}}^{(2)}$, we obtain
$$
c_{\bar{1}}^{+(2)}(\partial+\lambda,\mu)-(2\D+\L)\lambda(f_{\bar{1}}^{+(2)}(\partial+\lambda,\mu)-\delta^*_\mu)+\delta_{k+p,0}\b=0.
$$
Taking $\lambda=0$, we see that $c_{\bar{1}}^{+(2)}(\partial,\mu)=-\delta_{k+p,0}\b$.
As above, by considering the actions of operators $[J_0\,{}_\lambda\, G_k^+]\,{}_{\lambda+\mu}=G_k^+{}_{\lambda+\mu}$
and $[L_0\,{}_\lambda\, G_k^+]\,{}_{\lambda+\mu}=((k+\frac{1}{2}p)\lambda-p\mu)G_k^+{}_{\lambda+\mu}$ on $v_{\bar{1}}^{(2)}$,
one can show that $f_{\bar{1}}^{+(2)}(\partial,\lambda)=0$.
Hence, the action of $G_k^+$ on $v_{\bar{1}}^{(2)}$ has the required form.
\QED

\begin{lemm}\label{action-J-on-M-1}
Let $\gamma$ be as in Lemma~\ref{action-G-on-M}.
We have $\gamma=0$, and $J_k\,{}_\lambda\, v_{\bar{0}}^{(1)}=J_k\,{}_\lambda\, v_{\bar{1}}^{(1)}=0$.
\end{lemm}

\ni{\it Proof.}\ \
We first show that $J_k\,{}_\lambda\, v_{\bar{1}}^{(1)}=0$.
Considering the action of $[J_k\,{}_\lambda\, G_0^{+}]\,{}_{\lambda+\mu}=G_k^{+}\,{}_{\lambda+\mu}$ on $v_{\bar{0}}^{(1)}$,
by Lemma~\ref{action-G-on-M}(1), we obtain $e_{\bar{1}}^{(1)}(\partial,\lambda)=0$ and
\begin{equation}\label{JkG0-V01}
b_{\bar{1}}^{(1)}(\partial,\lambda)-b_{\bar{0}}^{(1)}(\partial+\mu,\lambda)+(2\D+\L)\lambda e_{\bar{0}}^{(1)}(\partial+\mu,\lambda)=\frac{1}{\sqrt{p}}\delta_{\lambda+\mu}^*.
\end{equation}
Then, considering the action of $[J_k\,{}_\lambda\, J_k]\,{}_{\lambda+\mu}=0$ on $v_{\bar{1}}^{(1)}$, we obtain
$$
b_{\bar{1}}^{(1)}(\partial,\lambda)b_{\bar{1}}^{(1)}(\partial+\lambda,\mu) = b_{\bar{1}}^{(1)}(\partial,\mu)b_{\bar{1}}^{(1)}(\partial+\mu,\lambda).
$$
Comparing the coefficients of $\lambda$ in the above equations,
we see that $b_{\bar{1}}^{(1)}(\partial,\lambda)$ is independent of the variable $\partial$,
and so we can denote $b_{\bar{1}}^{(1)}(\lambda)=b_{\bar{1}}^{(1)}(\partial,\lambda)$.
Furthermore, applying $[L_0\,{}_\lambda\, J_k]\,{}_{\lambda+\mu}=(k\lambda-p\mu)J_k{}_{\lambda+\mu}$ on $v_{\bar{1}}^{(1)}$, we obtain
$-p\mu b_{\bar{1}}^{(1)}(\mu)=(k\lambda-p\mu) b_{\bar{1}}^{(1)}(\lambda+\mu)$.
Taking $\mu=0$, we see that $b_{\bar{1}}^{(1)}(\partial,\lambda)=b_{\bar{1}}^{(1)}(\lambda)=0$. Hence, $J_k\,{}_\lambda\, v_{\bar{1}}^{(1)}=0$.

Now, taking $\mu=0$ and $\lambda=0$ in \eqref{JkG0-V01}, respectively, we obtain
\begin{eqnarray}
\label{JkG0-V01-mu=0} b_{\bar{0}}^{(1)}(\partial,\lambda) &\!\!\!=\!\!\!& (2\D+\L)\lambda e_{\bar{0}}^{(1)}(\partial,\lambda)-\frac{1}{\sqrt{p}}\delta_{\lambda}^*,\\
\label{JkG0-V01-lambda=0} b_{\bar{0}}^{(1)}(\partial+\mu, 0) &\!\!\!=\!\!\!& -\frac{1}{\sqrt{p}}\delta_{\mu}^*.
\end{eqnarray}
Considering the action of $[J_0\,{}_\lambda\, J_k]\,{}_{\lambda+\mu}=0$ on $v_{\bar{0}}^{(1)}$
and equating the coefficients of $v_{\bar{0}}^{(1)}$, we obtain
\begin{eqnarray}
\label{J0Jk-V01b} \L b_{\bar{0}}^{(1)}(\partial+\lambda,\mu)- \L b_{\bar{0}}^{(1)}(\partial,\mu)+(2\D+\L)\lambda e_{\bar{0}}^{(1)}(\partial+\lambda,\mu) &\!\!\!=\!\!\!& 0.
\end{eqnarray}

Next, we show that $\gamma=0$, or equivalently $\delta_{\lambda}^*=0$. Note first that from the definition of $\delta_{\lambda}^*$,
we only need to consider the case $k=p$ and $\D=\frac{3}{2}\L+1$.
By \eqref{JkG0-V01-lambda=0}, we see that $\delta_{\mu}^*$ is a polynomial on $\partial+\mu$.
Recall again the definition of $\delta_{\mu}^*$, we may further assume that $\L=0$, and so $\D=1$.
Then, by \eqref{J0Jk-V01b}, we have $e_{\bar{0}}^{(1)}(\partial,\mu)=0$.
On one hand, this, together with \eqref{JkG0-V01-mu=0}, implies that
$b_{\bar{0}}^{(1)}(\partial,\lambda)=-\frac{1}{\sqrt{p}}\delta_{\lambda}^*=-\frac{1}{\sqrt{p}}\gamma(\partial+\lambda+\a)$.
On the other hand, applying $[J_k\,{}_\lambda\, J_k]\,{}_{\lambda+\mu}=0$ on $v_{\bar{0}}^{(1)}$, as above,
we must have that $b_{\bar{0}}^{(1)}(\partial,\lambda)$ is independent of the variable $\partial$. Hence, $\gamma=0$.

Applying $[L_0\,{}_\lambda\, J_k]\,{}_{\lambda+\mu}=(k\lambda-p\mu)J_k{}_{\lambda+\mu}$ on $v_{\bar{0}}^{(1)}$, and then taking $\mu=0$,
we obtain
\begin{eqnarray}
\label{L0Jk-V01-b} p(\partial+\D\lambda+\a)\frac{b_{\bar{0}}^{(1)}(\partial+\lambda,0)-b_{\bar{0}}^{(1)}(\partial,0)}{\lambda}
+p(\D+\frac{\L}{2})\lambda e_{\bar{0}}^{(1)}(\partial+\lambda,0) &\!\!\!=\!\!\!& k b_{\bar{0}}^{(1)}(\partial,\lambda),\\
\label{L0Jk-V01-e} p(\partial+\D\lambda+\a)\frac{e_{\bar{0}}^{(1)}(\partial+\lambda,0)-e_{\bar{0}}^{(1)}(\partial,0)}{\lambda}
+p e_{\bar{0}}^{(1)}(\partial+\lambda,0) &\!\!\!=\!\!\!& k e_{\bar{0}}^{(1)}(\partial,\lambda).
\end{eqnarray}
Taking $\lambda\rightarrow0$ in \eqref{L0Jk-V01-e}, we obtain
\begin{equation}\label{L0Jk-V01-e-lambda=0}
p(\partial+\a)\frac{d}{d\partial}e_{\bar{0}}^{(1)}(\partial,0)=(k-p)e_{\bar{0}}^{(1)}(\partial,0).
\end{equation}

At last, we show that $J_k\,{}_\lambda\, v_{\bar{0}}^{(1)}=0$ in two cases.
First, we consider the case $2\D+\L\ne0$.
On one hand, by \eqref{JkG0-V01-mu=0} (recall that $\delta_{\lambda}^*=0$), we have $b_{\bar{0}}^{(1)}(\partial,\lambda) = (2\D+\L)\lambda e_{\bar{0}}^{(1)}(\partial,\lambda)$.
In particular, $b_{\bar{0}}^{(1)}(\partial,0)=0$.
On the other hand, by \eqref{L0Jk-V01-b}, we have $b_{\bar{0}}^{(1)}(\partial,\lambda)=\frac{p}{k}(\D+\frac{\L}{2})\lambda e_{\bar{0}}^{(1)}(\partial+\lambda,0)$.
Hence, we have $\frac{p}{2k}e_{\bar{0}}^{(1)}(\partial+\lambda,0)=e_{\bar{0}}^{(1)}(\partial,\lambda)$.
If $p\ne 2k$, by comparing the coefficients of $\partial$ on both sides, we must have that $e_{\bar{0}}^{(1)}(\partial,\lambda)=0$,
and thus $b_{\bar{0}}^{(1)}(\partial,\lambda)=0$. Namely, $J_k\,{}_\lambda\, v_{\bar{0}}^{(1)}=0$.
If $p= 2k$, by \eqref{L0Jk-V01-e-lambda=0}, we must have $e_{\bar{0}}^{(1)}(\partial,0)=0$.
Then, by \eqref{L0Jk-V01-e} and then \eqref{JkG0-V01-mu=0}, we also have that $J_k\,{}_\lambda\, v_{\bar{0}}^{(1)}=0$.

Next, we suppose $2\D+\L=0$. First, by \eqref{JkG0-V01-mu=0}, we have $b_{\bar{0}}^{(1)}(\partial,\lambda) = 0$.
So, we only need to show $e_{\bar{0}}^{(1)}(\partial,\lambda) = 0$.
If $k=-p$, by \eqref{L0Jk-V01-e-lambda=0}, we must have $e_{\bar{0}}^{(1)}(\partial,0)=0$.
As above, by \eqref{L0Jk-V01-e} and then \eqref{JkG0-V01-mu=0}, we have $J_k\,{}_\lambda\, v_{\bar{0}}^{(1)}=0$.
If $k\ne-p$, then by Lemma~\ref{action-G-on-M} and $\delta_{\lambda}^*=0$, we see that the action of $G_k^+$ on $M$ is trivial.
Applying the operator $[G_k^{+}\,{}_\lambda\, G_0^{-}]\,{}_{\lambda+\mu}=(p\lambda-(2k+p)\mu)J_k\,{}_{\lambda+\mu}+2 L_k\,{}_{\lambda+\mu}$
on $v_{\bar{0}}^{(1)}$, we see that $L_k\,{}_{\lambda+\mu} v_{\bar{0}}^{(1)}=-\frac{1}{2}(p\lambda-(2k+p)\mu)J_k\,{}_{\lambda+\mu}v_{\bar{0}}^{(1)}$.
In particular, taking $\mu=0$, we have $L_k\,{}_{\lambda} v_{\bar{0}}^{(1)}=-\frac{1}{2}p\lambda J_k\,{}_{\lambda}v_{\bar{0}}^{(1)}$.
Then, applying $[L_k\,{}_\lambda\, J_0]\,{}_{\lambda+\mu}=-(k+p)\mu J_k\,{}_{\lambda+\mu}$ on $v_{\bar{0}}^{(1)}$,
we obtain
$$
\frac{p}{2}\L\lambda(e_{\bar{0}}^{(1)}(\partial+\mu,\lambda)-e_{\bar{0}}^{(1)}(\partial,\lambda))=-(k+p)\mu e_{\bar{0}}^{(1)}(\partial,\lambda+\mu).
$$
Taking $\lambda=0$, we see that $e_{\bar{0}}^{(1)}(\partial,\mu)=0$. This completes the proof.
\QED

\begin{lemm}\label{action-J-on-M-2}
Let $\b$ be as in Lemma~\ref{action-G-on-M}.
We have $J_k\,{}_\lambda\, v_{\bar{1}}^{(2)}=0$,  $J_k\,{}_\lambda\, v_{\bar{0}}^{(2)}=-\delta_{k+p,0}\frac{\b }{\sqrt{p}}v_{\bar{0}}^{(1)}$.
\end{lemm}

\ni{\it Proof.}\ \
We first show that $J_k\,{}_\lambda\, v_{\bar{0}}^{(2)}=-\delta_{k+p,0}\frac{\b }{\sqrt{p}}v_{\bar{0}}^{(1)}$.
Considering the action of $[J_k\,{}_\lambda\, G_0^{+}]\,{}_{\lambda+\mu}=G_k^{+}\,{}_{\lambda+\mu}$ on $v_{\bar{0}}^{(2)}$,
by Lemma~\ref{action-G-on-M}(2) and Lemma~\ref{action-J-on-M-1}, we obtain
$$
b_{\bar{0}}^{(2)}(\partial+\mu,\lambda)+(2\D+\L)\lambda e_{\bar{0}}^{(2)}(\partial+\mu,\lambda)=-\delta_{k+p,0}\frac{\b}{\sqrt{p}}.
$$
Taking $\lambda=\mu=0$, we have $b_{\bar{0}}^{(2)}(\partial,0)=-\delta_{k+p,0}\frac{\b}{\sqrt{p}}$.
Considering the action of $[J_0\,{}_\lambda\, J_k]\,{}_{\lambda+\mu}=0$ on $v_{\bar{0}}^{(2)}$
and equating the coefficients of $v_{\bar{0}}^{(1)}$, by Lemma~\ref{action-J-on-M-1}, we obtain
\begin{equation}\label{J0Jk-V01b}
\L b_{\bar{0}}^{(2)}(\partial+\lambda,\mu)- \L b_{\bar{0}}^{(2)}(\partial,\mu) +(2\D+\L)\lambda e_{\bar{0}}^{(1)}(\partial+\lambda,\mu)=0.
\end{equation}
Applying $[L_0\,{}_\lambda\, J_k]\,{}_{\lambda+\mu}=(k\lambda-p\mu)J_k{}_{\lambda+\mu}$ on $v_{\bar{0}}^{(2)}$, and then taking $\mu=0$,
we obtain (recall that we have shown $b_{\bar{0}}^{(2)}(\partial,0)=-\delta_{k+p,0}\frac{\b}{\sqrt{p}}$)
\begin{eqnarray}
\label{L0Jk-V02-b} -p(\D+\frac{\L}{2})\lambda e_{\bar{0}}^{(2)}(\partial+\lambda,0)+ \delta_{k+p,0}\sqrt{p}\b
 &\!\!\!=\!\!\!& k b_{\bar{0}}^{(2)}(\partial,\lambda),\\
\label{L0Jk-V02-e} p(\partial+(\D+1)\lambda+\a)\frac{e_{\bar{0}}^{(2)}(\partial+\lambda,0)-e_{\bar{0}}^{(2)}(\partial,0)}{\lambda}
 &\!\!\!=\!\!\!& k e_{\bar{0}}^{(2)}(\partial,\lambda).
\end{eqnarray}
Taking $\lambda\rightarrow0$ in \eqref{L0Jk-V02-e}, we obtain
\begin{equation}\label{L0Jk-V02-e-lambda=0}
p(\partial+\a)\frac{d}{d\partial}e_{\bar{0}}^{(2)}(\partial,0)=k e_{\bar{0}}^{(2)}(\partial,0).
\end{equation}
By \eqref{J0Jk-V01b} with $\mu=0$, we have $(2\D+\L)\lambda e_{\bar{0}}^{(1)}(\partial+\lambda,0)=0$.
This, together with \eqref{L0Jk-V02-b}, implies that
$b_{\bar{0}}^{(2)}(\partial,\lambda)=\delta_{k+p,0}\frac{\sqrt{p}\b}{k}=-\delta_{k+p,0}\frac{\b}{\sqrt{p}}$.
Let $\deg\, e_{\bar{0}}^{(2)}(\partial,0)=n$.
Assume that $n\ge 1$. By \eqref{L0Jk-V02-e-lambda=0}, we must have $k=pn$, which implies $k+p=p(n+1)\ne0$.
Hence, $b_{\bar{0}}^{(2)}(\partial,\lambda)=0$.
Then, applying $[J_k\,{}_\lambda\, J_k]\,{}_{\lambda+\mu}=0$ on $v_{\bar{0}}^{(2)}$, as in Lemma~\ref{action-J-on-M-1},
we must have that $e_{\bar{0}}^{(2)}(\partial,\lambda)$ is independent of the variable $\partial$.
In particular, $e_{\bar{0}}^{(2)}(\partial,0)$ is independent of the variable $\partial$, which contradicts to the assumption $\deg\, e_{\bar{0}}^{(2)}(\partial,0)=n\ge1$.
Thus, $n=0$. By \eqref{L0Jk-V02-e-lambda=0}, we have $e_{\bar{0}}^{(2)}(\partial,0)=0$.
Furthermore, by \eqref{L0Jk-V02-e}, we have $e_{\bar{0}}^{(2)}(\partial,\lambda)=0$.
Hence, the action of $J_k$ on $v_{\bar{0}}^{(2)}$ has the required form.

Next, we show that $J_k\,{}_\lambda\, v_{\bar{1}}^{(2)}=0$.
Considering the action of $[J_k\,{}_\lambda\, G_0^{+}]\,{}_{\lambda+\mu}=G_k^{+}\,{}_{\lambda+\mu}$ on $v_{\bar{1}}^{(2)}$, and then taking $\mu=0$,
by Lemma~\ref{action-G-on-M}(2), Lemma~\ref{action-J-on-M-1} and $J_k\,{}_\lambda\, v_{\bar{0}}^{(2)}=-\delta_{k+p,0}\frac{\b }{\sqrt{p}}v_{\bar{0}}^{(1)}$,
one can first show that $e_{\bar{1}}^{(2)}(\partial,\lambda)=0$.
Considering  further the action of $[J_0\,{}_\lambda\, J_k]\,{}_{\lambda+\mu}=0$ on $v_{\bar{1}}^{(2)}$, and then taking $\lambda=0$,
we obtain $b_{\bar{1}}^{(2)}(\partial,\mu)=0$. Hence, $J_k\,{}_\lambda\, v_{\bar{1}}^{(2)}=0$.
\QED
\vskip15pt

\ni{\it Continuation of the proof of Theorem \ref{thm-rank-2+2}.}\ \
If $k\ne -p$, then $k+p\ne 0$. By Lemmas~\ref{action-G-on-M}--\ref{action-J-on-M-2}, we see that the actions of $G_k^+$ and $J_k$ on $M$ are trivial.
Then, by relations $[J_k\,{}_\lambda\, G_0^{-}]=- G_k^{-}$ and $[G_0^{+}\,{}_\lambda\, G_k^{-}]=(p\partial+2(k+p)\lambda)J_k+2 L_k$,
we see that the actions of $G_k^{-}$ and $L_k$ on $M$ are also trivial.
This contradicts to the assumption that the action of $\KB_k$ is nontrivial.

Next, we assume that $k=-p$ (and thus $p$ is a negative integer).
Considering the action of  $[J_k\,{}_\lambda\, G_0^{-}]\,{}_{\lambda+\mu}=- G_k^{-}\,{}_{\lambda+\mu}$ on $M$,
by Lemmas~\ref{action-J-on-M-1} and \ref{action-J-on-M-2},  we obtain
\begin{equation}\label{action-G-minus-on-M}
G_k^{-}\,{}_\lambda\, v_{\bar{0}}^{(1)}=G_k^{-}\,{}_\lambda\, v_{\bar{1}}^{(2)}=0, \quad
G_k^{-}\,{}_\lambda\, v_{\bar{0}}^{(2)}=-\b v_{\bar{1}}^{(2)}, \quad G_k^{-}\,{}_\lambda\, v_{\bar{1}}^{(1)}=-\b v_{\bar{0}}^{(1)}.
\end{equation}
Then, applying $[G_0^{+}\,{}_\lambda\, G_k^{-}]\,{}_{\lambda+\mu}=-p(\lambda+\mu) J_k\,{}_{\lambda+\mu}+2 L_k\,{}_{\lambda+\mu}$ on $M$,
by Lemmas~\ref{action-J-on-M-1} and \ref{action-J-on-M-2},  we obtain
\begin{equation}\label{action-L-on-M}
L_k\,{}_\lambda\, v_{s}^{\ell}=-\frac{\sqrt{p}}{2}\b (v_{s}^{\ell}+\delta_{s,\bar{0}}\delta_{\ell,2}\lambda v_{\bar{0}}^{(1)}), \ \ s\in\Z/2\Z,\, \ell=1,2.
\end{equation}

If $p\le -2$, then $k\ge 2$.
Following the arguments for case $k+p\ne0$, one can show that the action of $\KB_{k-1}$ on $M$ is trivial.
Then, by relations
\begin{eqnarray*}
[L_1\, {}_\lambda \, L_{k-1}] &\!\!\!=\!\!\!& ((1+p)\partial+p\lambda) L_{k}, \\[-2pt]
[L_1\, {}_\lambda \, J_{k-1}] &\!\!\!=\!\!\!& (1+p)\partial J_{k}, \\[-2pt]
[L_1\, {}_\lambda \, G^{\pm}_{k-1}] &\!\!\!=\!\!\!& ((1+p)\partial+\frac{1}{2}p\lambda) G^{\pm}_{k},
\end{eqnarray*}
we see that the action of $\KB_{k}$ on $M$ is also trivial, a contradiction.
Hence, $p=-1$, and by Lemmas~\ref{action-G-on-M}--\ref{action-J-on-M-2}, \eqref{action-G-minus-on-M} and \eqref{action-L-on-M},
we have \eqref{KB-non-trivial-extension}. This completes the proof.
\QED

\subsection{Composition factors}

Let $M$ be a conformal module over $\KB$. In this subsection, we use $\widetilde{M}$ to denote the same module with reversed parity.
Next, we determine the composition factors of all the free conformal $\KB$-modules of small rank obtained
in Theorems~\ref{thm-rank-1+1} and \ref{thm-rank-2+2}.
In fact, for rank $(1+1)$ $\KB$-modules in Theorem~\ref{thm-rank-1+1}, one can easily give their simplicities and composition factors.

\begin{prop}\label{prop-simplicity-1+1}
Let $M$ be a conformal $\KB$-module in Theorem~\ref{thm-rank-1+1}.
\baselineskip1pt\lineskip7pt\parskip-1pt
\begin{itemize}\parskip-1pt
  \item[{\rm(1)}] If $M\cong V_{\D,\a}^{(1)}$, then $M$ is simple if and only if $\D\ne 0$. Furthermore, $V_{0,\a}^{(1)}$
  contains a unique nontrivial submodule $\C[\partial](\partial+\a)v_{\bar{0}}\oplus\C[\partial](\frac{1}{2}v_{\bar{1}})\cong \widetilde{V_{\frac{1}{2},\a}^{(2)}}$,
  and the quotient is an even trivial module $\C c_{-\a}$.
  \item[{\rm(2)}] If $M\cong V_{\D,\a}^{(2)}$, then $M$ is simple if and only if $\D\ne 0$. Furthermore, $V_{0,\a}^{(2)}$
  contains a unique nontrivial submodule $\C[\partial](\partial+\a)v_{\bar{0}}\oplus\C[\partial](\frac{1}{2}v_{\bar{1}})\cong \widetilde{V_{\frac{1}{2},\a}^{(1)}}$,
  and the quotient is an even trivial module $\C c_{-\a}$.
\end{itemize}
In addition, the above conclusions still hold if we reverse the parity of all $\KB$-modules.
\end{prop}

Nevertheless, for rank $(2+2)$ $\KB$-modules in Theorem~\ref{thm-rank-2+2}, the situation is not so obvious.

\begin{prop}\label{prop-simplicity-2+2}
Let $M$ be a conformal $\KB$-module in Theorem~\ref{thm-rank-2+2}.
\baselineskip1pt\lineskip7pt\parskip-1pt
\begin{itemize}\parskip-1pt
  \item[{\rm(1)}] If $p\ne-1$ and $M\cong V_{\D,\L,\a}$, then $M$ is simple if and only if $2\D\pm\L\ne 0$. Furthermore,
  \begin{itemize}\parskip-1pt
  \item[\rm{(i)}] $V_{\D,-2\D,\a}$ contains a submodule isomorphic to $\widetilde{V_{\D+\frac{1}{2},\a}^{(1)}}$,
  and the quotient is $V_{\D,\a}^{(1)}$;
  \item[\rm{(ii)}] $V_{\D,2\D,\a}$ contains a submodule isomorphic to $\widetilde{V_{\D+\frac{1}{2},\a}^{(2)}}$,
  and the quotient is $V_{\D,\a}^{(2)}$.
  \end{itemize}
  \item[{\rm(2)}] If $p=-1$ and $M\cong V_{\D,\L,\a,\b}$, then $M$ is simple if and only if $(2\D\pm\L,\b)\ne (0,0)$. Furthermore,
  \begin{itemize}\parskip-1pt
  \item[\rm{(i)}] $V_{\D,-2\D,\a,0}$ contains a submodule isomorphic to $\widetilde{V_{\D+\frac{1}{2},\a}^{(1)}}$,
  and the quotient is $V_{\D,\a}^{(1)}$;
  \item[\rm{(ii)}] $V_{\D,2\D,\a,0}$ contains a submodule isomorphic to $\widetilde{V_{\D+\frac{1}{2},\a}^{(2)}}$,
  and the quotient is $V_{\D,\a}^{(2)}$.
  \end{itemize}
\end{itemize}
In addition, the above conclusions still hold if we reverse the parity of all $\KB$-modules.
All composition factors of $M$ have multiplicity one.
The composition factors of $V_{\D,\pm 2\D,\a}$ and $V_{\D,\pm 2\D,\a,0}$ are listed as follows:
\begin{table}[h]
\centering\small
\subtable[Composition factors for cases (i)'s]{
\begin{tabular}{l|l}
\hline
 $V_{\D,-2\D,\a}$/$V_{\D,-2\D,\a,0}$  & Composition factors  \\
\hline\\[-12pt]
$\D\ne 0, -\frac{1}{2}$  & $V^{(1)}_{\D,\a}$, $\widetilde{V^{(1)}_{\D+\frac{1}{2},\a}}$  \\[5pt]
\hline\\[-12pt]
$\D= 0$ & $\widetilde{V^{(1)}_{\frac{1}{2},\a}}$, $\widetilde{V^{(2)}_{\frac{1}{2},\a}}$, $\C c_{-\a}$ \\[5pt]
\hline\\[-12pt]
$\D= -\frac{1}{2}$ & $V^{(1)}_{-\frac{1}{2},\a}$, $V^{(2)}_{\frac{1}{2},\a}$, $\widetilde{\C c_{-\a}}$ \\[5pt]
\hline
\end{tabular}}
\quad
\subtable[Composition factors for cases (ii)'s]{
\begin{tabular}{l|l}
\hline
 $V_{\D,2\D,\a}$/$V_{\D,2\D,\a,0}$  & Composition factors  \\
\hline\\[-12pt]
$\D\ne 0, -\frac{1}{2}$  & $V^{(2)}_{\D,\a}$, $\widetilde{V^{(2)}_{\D+\frac{1}{2},\a}}$  \\[5pt]
\hline\\[-12pt]
$\D= 0$ & $\widetilde{V^{(1)}_{\frac{1}{2},\a}}$, $\widetilde{V^{(2)}_{\frac{1}{2},\a}}$, $\C c_{-\a}$ \\[5pt]
\hline\\[-12pt]
$\D= -\frac{1}{2}$ & $V^{(1)}_{\frac{1}{2},\a}$, $V^{(2)}_{-\frac{1}{2},\a}$, $\widetilde{\C c_{-\a}}$ \\[5pt]
\hline
\end{tabular}}
\end{table}
\end{prop}

\ni{\it Proof.}\ \
(1) First, recall that $V_{\D,\L,\a}$ is a trivial extension of the $K_2$-module $K_{\D,\L,\a}$.
By Remark~\ref{remark-K2}, $M$ is simple if and only if $2\D\pm\L\ne 0$.

If $2\D+\L=0$, then $M\cong V_{\D,-2\D,\a}$.
Let $v_{\bar{0}}=v_{\bar{0}}^{(2)}, v_{\bar{1}}=v_{\bar{1}}^{(2)}$. By \eqref{K2-module-rk2+2} with $\L=-2\D$,
it is easy to observe that $M_1=\C[\partial]v_{\bar{0}}\oplus\C[\partial]v_{\bar{1}}$ is a submodule isomorphic to $\widetilde{V_{\D+\frac{1}{2},\a}^{(1)}}$,
and $M/M_1\cong V_{\D,\a}^{(1)}$.
If $\D\ne 0, -\frac{1}{2}$, then both $M_1$ and $M/M_1$ are irreducible.
If $\D=0$, then $M_1\cong \widetilde{V^{(1)}_{\frac{1}{2},\a}}$ is irreducible,
while $M/M_1\cong V_{0,\a}^{(1)}$ has composition factors $\widetilde{V^{(2)}_{\frac{1}{2},\a}}$ and $\C c_{-\a}$ by Proposition~\ref{prop-simplicity-1+1}(1).
If $\D=-\frac{1}{2}$, then $M/M_1\cong V_{-\frac{1}{2},\a}^{(1)}$ is irreducible,
while $M_1\cong \widetilde{V^{(1)}_{0,\a}}$ has composition factors $V^{(2)}_{\frac{1}{2},\a}$ and $\widetilde{\C c_{-\a}}$ by
the parity reverse version of Proposition~\ref{prop-simplicity-1+1}(1).

If $2\D-\L=0$, then $M\cong V_{\D,2\D,\a}$. Although any trivial combination of $v_{\bar{0}}^{(1)},v_{\bar{0}}^{(2)},v_{\bar{1}}^{(1)},v_{\bar{1}}^{(2)}$
is not a rank $(1+1)$ submodule of $\KB$, we have a feeling that there should be certain symmetry between cases $2\D-\L=0$ and $2\D+\L=0$.
So, let us consider the structure of $\widetilde{V_{\D+\frac{1}{2},\a}^{(2)}}$ (cf.~\eqref{K2-module-rk1+1-2}):
\begin{equation*}
\left\{\begin{array}{ll}
L_0\,{}_\lambda\, v_{\bar{1}}=p(\partial+(\D+\frac{1}{2}) \lambda+\a)v_{\bar{1}}, &  L_0\,{}_\lambda\, v_{\bar{0}}=p(\partial+(\D+1) \lambda+\a)v_{\bar{0}},  \\[4pt]
J_0\,{}_\lambda\, v_{\bar{1}}=(2 \D+ 1) v_{\bar{1}}, &  J_0\,{}_\lambda\, v_{\bar{0}}=2\D v_{\bar{0}},  \\[4pt]
G_0^+\,{}_\lambda\, v_{\bar{1}}=0, & G_0^+\,{}_\lambda\, v_{\bar{0}}=2\sqrt{p}(\partial+(2\D+1)\lambda+\a)v_{\bar{1}}, \\[4pt]
G_0^-\,{}_\lambda\, v_{\bar{1}}=\sqrt{p}v_{\bar{0}}, & G_0^-\,{}_\lambda\, v_{\bar{0}}=0, \\[4pt]
X_i\,{}_\lambda\, v_{s}=0,
\end{array}\right.
\end{equation*}
where $X=L,J,G^\pm,\, i\ge 1,\, s\in\Z/2\Z$.
Comparing the above with \eqref{K2-module-rk2+2} with $\L=2\D$ (especially the action of $G_0^-$),
we naturally set $v_{\bar{0}}=2(\partial+\a)v_{\bar{0}}^{(1)}-v_{\bar{0}}^{(2)}, v_{\bar{1}}=v_{\bar{1}}^{(1)}$.
Then, one can check that $M_2=\C[\partial]v_{\bar{0}}\oplus\C[\partial]v_{\bar{1}}$ is a submodule isomorphic to $\widetilde{V_{\D+\frac{1}{2},\a}^{(2)}}$,
and $M/M_2\cong V_{\D,\a}^{(2)}$.
As in case $2\D+\L=0$, the composition factors of $M$ can be determined by considering the simplicities of $M_2$ and $M/M_2$
by Proposition~\ref{prop-simplicity-1+1}(2) or its parity reverse version.

(2) If $\b=0$, then the arguments are the same as those in (1). We only need to note that if $\b\ne0$, then $M$ is irreducible.
If this is not true, we may assume that $M'$ is a submodule of $M$.
If $M'$ has rank $(2+2)$, then $M/M'$ is a trivial module.
One can easily derive a contradiction from \eqref{KB-non-trivial-extension}.
If $M'$ has rank $(1+1)$, then $M/M'$ also has rank $(1+1)$. 
This is also impossible, since the action of $\KB_1$ on any rank $(1+1)$ module
is trivial by Theorem~\ref{thm-rank-1+1}.
\QED

\begin{rema}\label{remark-K2-CF}
\rm
Recall that $K_2\cong\KB_{[0]}$ (cf.~\eqref{equ-quotiont}) can be viewed as a quotient algebra of $\KB$.
Hence, the classification of composition factors of reducible $\KB$-modules $V_{0,\a}^{(1)}$, $V_{0,\a}^{(2)}$, $V_{\D,\pm 2\D,\a}$
in Proposition~\ref{prop-simplicity-1+1} and Proposition~\ref{prop-simplicity-2+2}(1),
is also true for  reducible $K_2$-modules $K_{0,\a}^{(1)}$, $K_{0,\a}^{(2)}$, $K_{\D,\pm 2\D,\a}$
(one only need to replace the symbols $V$'s by $K$'s in all the involved modules).
\end{rema}

\section{\large{Classification theorems}}

\subsection{Main result}

Our main result is the following.

\begin{theo}\label{thm-classification}
Let $M$ be a nontrivial FICM over $\KB$.
\baselineskip1pt\lineskip7pt\parskip-1pt
\begin{itemize}\parskip-1pt
  \item[{\rm(1)}] If $p\ne -1$, then,  up to parity change,  $M$ is isomorphic to $V^{(1)}_{\D,\a}$ or $V^{(2)}_{\D,\a}$ for some $\D,\a\in\C$ with $\D\ne 0$,
  or $V_{\D,\L,\a}$ for some $\D,\L,\a\in\C$ with $2\D\pm\L\ne 0$.
  \item[{\rm(2)}] If $p= -1$, then,  up to parity change,  $M$ is isomorphic to $V^{(1)}_{\D,\a}$ or $V^{(2)}_{\D,\a}$ for some $\D,\a\in\C$ with $\D\ne 0$,
  or $V_{\D,\L,\a,\b}$ for some $\D,\L,\a,\b\in\C$ with $(2\D\pm\L,\b)\ne (0,0)$.
\end{itemize}
\end{theo}

Let $M$ be a nontrivial FICM over $\KB$.
Using previous results, we shall prove that $M$ must be free of rank $(1+1)$ or $(2+2)$.
Then the main result will follow from Theorems~\ref{thm-rank-1+1} and \ref{thm-rank-2+2} and Propositions~\ref{prop-simplicity-1+1} and \ref{prop-simplicity-2+2}.
The following technical result \cite{CK} will be also used.

\begin{lemm}\label{lemma-general}
Let $\cal{L}$ be a Lie superalgebra with a descending sequence of subspaces
${\cal{L}}\supset{\cal{L}}_{0}\supset{\cal{L}}_{1}\supset\ldots$ and an element $T$ satisfying $[T,{\cal L}_n] = {\cal L}_{n-1}$ for
$n\ge 1$. Let $V$ be an $\cal{L}$-module and let
$$
V_n = \{v\in V\,|\, {\cal{L}}_n v = 0\},  \ \ \ n\in\Z_+.
$$
Suppose that $V_n\ne 0$ for $n\gg 0$, and that the minimal $N\in\Z_+$ for which $V_N\ne 0$ is positive.
Then $\C[T]V_N=\C[T]\otimes_{\,\C} V_N$.
In particular, $V_N$ is finite-dimensional if $V$ is a finitely generated $\C[T]$-module.
\end{lemm}

\begin{lemm}\label{free-of-rank-1+1/2+2}
The conformal $\KB$-module $M$ must be free of rank $(1+1)$ or $(2+2)$.
\end{lemm}
\ni{\it Proof.}\ \
Based on technical results prepared in Sections~3 and 4 (especially Theorem~\ref{thm-1}),
one can safely generalize the arguments in our previous papers \cite{SXY,X} to the case here.
For completeness, we still write down the details.

Note first that, by Theorem~\ref{thm-2}, the $\lambda$-actions of $L_i$, $J_i$, $G_i^\pm$ on $M$ are trivial for $i\gg 0$.
Suppose that $k\in\Z_+$ is the largest integer such that the $\lambda$-action of $\KB_k$ (with $\C[\partial]$-basis
$\{L_k,\,J_k,\,G_k^\pm\}$) on $M$ is nontrivial.
Then $M$ is simply a nontrivial FICM over $\KB_{[k]}$, where $\KB_{[k]}$ is defined by \eqref{equ-quotiont}.
Furthermore, by Proposition~\ref{prop-observation}, $M$
can be viewed as a module over the Lie superalgebra
${\cal{L}}={\cal A}(\KB_{[k]})^e$ satisfying
\begin{equation}\label{equiverlent}
\bar{L}_{i,m} v = \bar{J}_{i,n} v =\bar{G}_{i,t}^\pm v = 0\ \ \ \mbox{for}\ \ v\in M,\ \ 0\le i\le k,\ \ 0\ll m,n\in\Z,\ \ \frac{1}{2}\ll t\in\frac{1}{2}+\Z.
\end{equation}
For $z\in\Z_+$, let
$$
{\cal{L}}_z=\sp_{\C}\{\bar{L}_{i,m},\,\bar{J}_{i,n},\,\bar{G}_{i,t}^\pm\in {\cal{L}}\,|\,0\le i\le k,\, z-1\le m\in\Z, \, z\le n\in\Z, \, z-\frac{1}{2}\le t\in\frac{1}{2}+\Z \}.
$$
Then ${\cal{L}}_{0}={\cal A}(\KB_{[k]})$
and ${\cal{L}}\supset{\cal{L}}_{0}\supset{\cal{L}}_{1}\supset\ldots$. By the definition of extended annihilation superalgebra,
we see that the element ${T}\in{\cal{L}}$ satisfies $[{T},{\cal L}_z] = {\cal L}_{z-1}$ for $z\ge 1$.
Let
$$
M_z = \{v\in M\,|\, {\cal{L}}_z v = 0\},  \ \ z\in\Z_+.
$$
By \eqref{equiverlent}, $M_z\ne 0$ for $z\gg 0$. Assume that $N\in\Z_+$ is the smallest integer such that $M_N\ne \emptyset$.

Similar to our previous results for $\CB$ \cite{SXY} and $\SB$ \cite{X}, the case $N=0$ is impossible.

Next, consider the case $N\ge 1$.
By the definition of extended annihilation superalgebra and the shift used in the proof of Lemma~\ref{annihilation-algebra}, we have
that  $T-\frac{1}{p} \bar{L}_{0,-1}\in\cal{L}$ is an even central element,
and so $T-\frac{1}{p} \bar{L}_{0,-1}$ acts on $M$ as a scalar.
Therefore, ${\cal{L}}_{0}$ acts irreducibly on $M$.
Furthermore,  since
$$
\bar{L}_{i,-1}=\frac{1}{p}[\bar{L}_{i,0},\bar{L}_{0,-1}],\quad
\bar{J}_{i,0}=\frac{1}{p}[\bar{J}_{i,1},\bar{L}_{0,-1}],\quad
\bar{G}_{i,-\frac{1}{2}}^\pm=\frac{1}{p}[\bar{G}_{i,\frac{1}{2}}^\pm,\bar{L}_{0,-1}],
$$
we see that
the action of ${\cal{L}}_{0}$ is determined by ${\cal{L}}_{1}$ and $\bar{L}_{0,-1}$ (or equivalently, determined by ${\cal{L}}_{1}$ and $T$).
Note that $M_N$ is ${\cal{L}}_1$-invariant.
By the irreducibility of $M$ and Lemma~\ref{lemma-general}, we see that $M=\C[T]\otimes_{\,\C} M_N$ and  $M_N$ is a nontrivial irreducible finite-dimensional ${\cal{L}}_1$-module.

If $N=1$, by definition we see that $M_1$ is in fact a trivial ${\cal{L}}_1$-module, a contradiction.

If $N\ge 2$, by definition we see that $M_N$ can be viewed as a ${\cal{L}}_1/{\cal{L}}_N$-module.
Note that  ${\cal{L}}_1/{\cal{L}}_N\cong {\frak g}(k,N-2)$.  By Theorem~\ref{thm-1}, we must have that
the dimension of $M_N$ is either $(1|1)$ or $(2|2)$. Equivalently, as a conformal $\KB$-module, $M$ is free of rank $(1+1)$ or $(2+2)$.
\QED

\subsection{Applications}

Recall the definition \eqref{equ-quotiont-special} of ${\frak k}(n)$, we have ${\frak k}(n)=\langle\bar{L}_i,\,\bar{J}_i,\,\bar{G}^\pm_i\,|\,0\le i\le n\rangle$
with the following nontrivial $\lambda$-brackets ($i+j\le n$):
\begin{eqnarray}
\nonumber[\bar{L}_i\, {}_\lambda \, \bar{L}_j] &\!\!\!=\!\!\!& ((i-n)\partial+(i+j-2n)\lambda) \bar{L}_{i+j}, \\[-2pt]
\nonumber[\bar{L}_i\, {}_\lambda \, \bar{J}_j] &\!\!\!=\!\!\!& ((i-n)\partial+(i+j-n)\lambda) \bar{J}_{i+j}, \\[-2pt]
\nonumber[\bar{L}_i\, {}_\lambda \, \bar{G}^{\pm}_j] &\!\!\!=\!\!\!& ((i-n)\partial+(i+j-\frac{3}{2}n)\lambda) \bar{G}^{\pm}_{i+j}, \\[-2pt]
\nonumber[\bar{J}_i\, {}_\lambda \, \bar{G}^{\pm}_j] &\!\!\!=\!\!\!& \pm \bar{G}^{\pm}_{i+j}, \\[-2pt]
\nonumber[\bar{G}^{+}_i\, {}_\lambda \, \bar{G}^{-}_j] &\!\!\!=\!\!\!& ((2i-n)\partial+2(i+j-n)\lambda) \bar{J}_{i+j}+2 \bar{L}_{i+j}.
\end{eqnarray}
Clearly, the following two $\C[\partial]$-modules are rank $(1+1)$ conformal modules over ${\frak k}(n)$
(here, we adopt the same notations as in (I-1) and (I-2) for $\KBB(-n)$, see Subsection~4.2).
\baselineskip1pt\lineskip7pt\parskip-1pt
\begin{itemize}\parskip-1pt
  \item[{\rm(i-1)}] $V_{\D,\a}^{(1)}=\C[\partial]v_{\bar{0}}\oplus\C[\partial]v_{\bar{1}}$ with
\begin{equation}\label{kn-trivial-extension-1}
\left\{\begin{array}{ll}
\bar{L}_0\,{}_\lambda\, v_{\bar{0}}=-n(\partial+\D \lambda+\a)v_{\bar{0}}, &  \bar{L}_0\,{}_\lambda\, v_{\bar{1}}=-n(\partial+(\D+\frac{1}{2}) \lambda+\a)v_{\bar{1}},  \\[4pt]
\bar{J}_0\,{}_\lambda\, v_{\bar{0}}=-2 \D v_{\bar{0}}, &  \bar{J}_0\,{}_\lambda\, v_{\bar{1}}=(1-2\D) v_{\bar{1}},  \\[4pt]
\bar{G}_0^+\,{}_\lambda\, v_{\bar{0}}=\sqrt{-n}v_{\bar{1}}, & \bar{G}_0^+\,{}_\lambda\, v_{\bar{1}}=0, \\[4pt]
\bar{G}_0^-\,{}_\lambda\, v_{\bar{0}}=0, & \bar{G}_0^-\,{}_\lambda\, v_{\bar{1}}=2\sqrt{-n}(\partial+2\D\lambda+\a)v_{\bar{0}}, \\[4pt]
\bar{X}_i\,{}_\lambda\, v_{s}=0,\, 1\le i\le n,\, s\in\Z/2\Z,
\end{array}\right.
\end{equation}
where $X=L,\,J,\,G^\pm$ and $\D,\,\a\in\C$;
  \item[{\rm(i-2)}] $V_{\D,\a}^{(2)}=\C[\partial]v_{\bar{0}}\oplus\C[\partial]v_{\bar{1}}$ with
\begin{equation}\label{kn-trivial-extension-2}
\left\{\begin{array}{ll}
\bar{L}_0\,{}_\lambda\, v_{\bar{0}}=-n(\partial+\D \lambda+\a)v_{\bar{0}}, &  \bar{L}_0\,{}_\lambda\, v_{\bar{1}}=-n(\partial+(\D+\frac{1}{2}) \lambda+\a)v_{\bar{1}},  \\[4pt]
\bar{J}_0\,{}_\lambda\, v_{\bar{0}}=2 \D v_{\bar{0}}, &  \bar{J}_0\,{}_\lambda\, v_{\bar{1}}=(2\D-1) v_{\bar{1}},  \\[4pt]
\bar{G}_0^+\,{}_\lambda\, v_{\bar{0}}=0, & \bar{G}_0^+\,{}_\lambda\, v_{\bar{1}}=2\sqrt{-n}(\partial+2\D\lambda+\a)v_{\bar{0}}, \\[4pt]
\bar{G}_0^-\,{}_\lambda\, v_{\bar{0}}=\sqrt{-n}v_{\bar{1}}, & \bar{G}_0^-\,{}_\lambda\, v_{\bar{1}}=0,\\[4pt]
\bar{X}_i\,{}_\lambda\, v_{s}=0,\, 1\le i\le n,\, s\in\Z/2\Z,
\end{array}\right.
\end{equation}
where $X=L,\,J,\,G^\pm$ and $\D,\,\a\in\C$;
\end{itemize}
The following $V_{\D,\L,\a}$ is a rank $(2+2)$ conformal module over ${\frak k}(n)$.
One can check that the following $V_{\D,\L,\a,\b}$ is a more general rank $(2+2)$ conformal module over ${\frak k}(1)$.
(here, we adopt the same notations as in (II-1) for $\KBB(-n)$ and (II-2) for $\KBB(-1)$, see Subsection~4.3.)
\begin{itemize}\parskip-1pt
  \item[{\rm(ii-1)}] $V_{\D,\L,\a}=\C[\partial]v_{\bar{0}}^{(1)}\oplus\C[\partial]v_{\bar{0}}^{(2)}\oplus\C[\partial]v_{\bar{1}}^{(1)}\oplus\C[\partial]v_{\bar{1}}^{(2)}$ with
\begin{equation}\label{kn-trivial-extension-3}
\left\{\begin{array}{ll}
\bar{L}_0\,{}_\lambda\, v_{\bar{0}}^{(1)}=-n(\partial+\D \lambda+\a)v_{\bar{0}}^{(1)},
\quad\ \   \bar{L}_0\,{}_\lambda\, v_{\bar{1}}^{(\ell)}=-n(\partial+(\D+\frac{1}{2}) \lambda+\a)v_{\bar{1}}^{(\ell)},\ \ell=1,2,  \\[4pt]
\bar{L}_0\,{}_\lambda\, v_{\bar{0}}^{(2)}=-n(\partial+(\D+1) \lambda+\a)v_{\bar{0}}^{(2)}-n(\D+\frac{\L}{2})\lambda^2 v_{\bar{0}}^{(1)},  \\[4pt]
\bar{J}_0\,{}_\lambda\, v_{\bar{0}}^{(1)}=\L v_{\bar{0}}^{(1)}, \qquad\qquad\qquad\qquad\   \bar{J}_0\,{}_\lambda\, v_{\bar{1}}^{(1)}=(\L+1) v_{\bar{1}}^{(1)},  \\[4pt]
\bar{J}_0\,{}_\lambda\, v_{\bar{0}}^{(2)}=\L v_{\bar{0}}^{(2)}+(2\D+\L)\lambda v_{\bar{0}}^{(1)}, \quad  \bar{J}_0\,{}_\lambda\, v_{\bar{1}}^{(2)}=(\L-1) v_{\bar{1}}^{(2)},  \\[4pt]
\bar{G}_0^+\,{}_\lambda\, v_{\bar{0}}^{(1)}=\sqrt{-n}v_{\bar{1}}^{(1)}, \qquad\qquad\qquad\ \bar{G}_0^+\,{}_\lambda\, v_{\bar{1}}^{(1)}=0, \\[4pt]
\bar{G}_0^+\,{}_\lambda\, v_{\bar{0}}^{(2)}=-\sqrt{-n}(2\D+\L)\lambda v_{\bar{1}}^{(1)},
\quad  \bar{G}_0^+\,{}_\lambda\, v_{\bar{1}}^{(2)}=\sqrt{-n}(2\D+\L)\lambda v_{\bar{0}}^{(1)} + \sqrt{-n} v_{\bar{0}}^{(2)}, \\[4pt]
\bar{G}_0^-\,{}_\lambda\, v_{\bar{0}}^{(1)}=\sqrt{-n}v_{\bar{1}}^{(2)},
\qquad\qquad\qquad\  \bar{G}_0^-{}_\lambda\, v_{\bar{1}}^{(1)}\!=\!\sqrt{-n}(2\partial\!+\!(2\D\!-\!\L)\lambda\!+\!2\a)v_{\bar{0}}^{(1)}\!\!-\!\!\sqrt{-n}v_{\bar{0}}^{(2)},\\[4pt]
\bar{G}_0^-\,{}_\lambda\, v_{\bar{0}}^{(2)}=\sqrt{-n}(2\partial\!+\!(2\D\!-\!\L\!+\!2)\lambda\!+\!2\a)v_{\bar{1}}^{(2)},
\ \, \bar{G}_0^-\,{}_\lambda\, v_{\bar{1}}^{(2)}=0,\\[4pt]
\bar{X}_i\,{}_\lambda\, v_{s}^{(\ell)}=0,\, 1\le i\le n,\, s\in\Z/2\Z,\,\ell=1,2,
\end{array}\right.
\end{equation}
where $X=L,\,J,\,G^\pm$ and $\D,\,\L,\,\a\in\C$;
  \item[{\rm(ii-2)}] $V_{\D,\L,\a,\b}=\C[\partial]v_{\bar{0}}^{(1)}\oplus\C[\partial]v_{\bar{0}}^{(2)}\oplus\C[\partial]v_{\bar{1}}^{(1)}\oplus\C[\partial]v_{\bar{1}}^{(2)}$
  with \eqref{kn-trivial-extension-3} and
\begin{equation}\label{kn-non-trivial-extension}
\left\{\begin{array}{ll}
   \bar{L}_1\,{}_\lambda\, v_{\bar{0}}^{(1)}=-\frac{\sqrt{-1}}{2}\b v_{\bar{0}}^{(1)},
&  \bar{L}_1\,{}_\lambda\, v_{\bar{1}}^{(1)}=-\frac{\sqrt{-1}}{2}\b v_{\bar{1}}^{(1)},  \\[4pt]
   \bar{L}_1\,{}_\lambda\, v_{\bar{0}}^{(2)}=-\frac{\sqrt{-1}}{2}\b (v_{\bar{0}}^{(2)}+\lambda v_{\bar{0}}^{(1)}),
&  \bar{L}_1\,{}_\lambda\, v_{\bar{1}}^{(2)}=-\frac{\sqrt{-1}}{2}\b v_{\bar{1}}^{(2)},  \\[4pt]
   \bar{J}_1\,{}_\lambda\, v_{\bar{0}}^{(1)}=0, &  \bar{J}_1\,{}_\lambda\, v_{\bar{1}}^{(1)}=0,  \\[4pt]
   \bar{J}_1\,{}_\lambda\, v_{\bar{0}}^{(2)}=\sqrt{-1}\b v_{\bar{0}}^{(1)}, &  \bar{J}_1\,{}_\lambda\, v_{\bar{1}}^{(2)}=0,  \\[4pt]
   \bar{G}_1^+\,{}_\lambda\, v_{\bar{0}}^{(1)}=0, & \bar{G}_1^+\,{}_\lambda\, v_{\bar{1}}^{(1)}=0, \\[4pt]
   \bar{G}_1^+\,{}_\lambda\, v_{\bar{0}}^{(2)}=\b v_{\bar{1}}^{(1)}, & \bar{G}_1^+\,{}_\lambda\, v_{\bar{1}}^{(2)}=-\b v_{\bar{0}}^{(1)}, \\[4pt]
   \bar{G}_1^-\,{}_\lambda\, v_{\bar{0}}^{(1)}=0, & \bar{G}_1^-\,{}_\lambda\, v_{\bar{1}}^{(1)}=-\b v_{\bar{0}}^{(1)}, \\[4pt]
   \bar{G}_1^-\,{}_\lambda\, v_{\bar{0}}^{(2)}=-\b v_{\bar{1}}^{(2)}, & \bar{G}_1^-\,{}_\lambda\, v_{\bar{1}}^{(2)}=0,
\end{array}\right.
\end{equation}
where $\D,\,\L,\,\a,\b\in\C$.
\end{itemize}

Since ${\frak k}(n)$ is a quotient algebra of $\KBB(-n)$ (cf.~\eqref{equ-sn} and \eqref{equ-quotiont-special}),
by Theorems~\ref{thm-rank-1+1} and \ref{thm-rank-2+2}, Propositions~\ref{prop-simplicity-1+1} and \ref{prop-simplicity-2+2}, we have that

\begin{coro}\label{classification-small-rank}
Let $M$ be a nontrivial free conformal module of rank $(1+1)$ or $(2+2)$ over ${\frak k}(n)$.
\baselineskip1pt\lineskip7pt\parskip-1pt
\begin{itemize}\parskip-1pt
  \item[{\rm(1)}] If $n>1$, then, up to parity change, $M$ is isomorphic to $V_{\D,\a}^{(1)}$ or $V_{\D,\a}^{(2)}$  for some $\D,\a\in\C$, or $V_{\D,\L,\a}$ for some $\D,\L,\a\in\C$.
  \item[{\rm(2)}] If $n=1$, then, up to parity change, $M$ is isomorphic to $V^{(1)}_{\D,\a}$ or $V^{(2)}_{\D,\a}$ for some $\D,\a\in\C$, or $V_{\D,\L,\a,\b}$ for some $\D,\L,\a,\b\in\C$.
\end{itemize}
Furthermore, for the above modules we have
the same assertions on simplicities and composition factors as those for $\KBB(-n)$-modules in Propositions~\ref{prop-simplicity-1+1} and \ref{prop-simplicity-2+2}.
\end{coro}

Furthermore,
by Theorem~\ref{thm-classification}, we have that

\begin{coro}\label{classification-kn}
The irreducible modules in Corollary~\ref{classification-small-rank}
exhaust all nontrivial finite irreducible conformal modules over ${\frak k}(n)$.
\end{coro}

\vskip10pt

\small \ni{\bf Acknowledgement}
This work was supported by the Fundamental Research Funds for the Central Universities, China (2019QNA34).

\end{document}